\documentclass[a4paper]{article}
\usepackage{fullpage}
\usepackage{amsmath}
\usepackage{amsthm}
\usepackage{amssymb}
\usepackage{amsbsy}
\usepackage{mathrsfs}
\usepackage{amstext,hyperref}  

\usepackage{chngcntr}
\usepackage{changepage}

\usepackage{array,longtable}
\usepackage{enumerate}

\usepackage{fullpage}

\newtheorem{theorem}{Theorem}[section] 
\newtheorem{proposition}[theorem]{Proposition} 
\newtheorem{definition}[theorem]{Definition} 
 
\newtheorem{corollary}[theorem]{Corollary} 
\newtheorem{claim}{Claim}[theorem] 
 
\newtheorem{question}[theorem]{Question} 

\newtheorem*{thm-main}{Theorem~\ref{thm-main}}

\hypersetup{
	colorlinks=true, 
	linkcolor=red,
	citecolor=blue,
}

\usepackage{graphicx}

\makeatletter
\DeclareRobustCommand{\ovl}[1]{%
	\mathpalette\do@cev{#1}%
}
\newcommand{\do@cev}[2]{%
	\fix@cev{#1}{+}%
	\reflectbox{$\m@th#1\ovr{\reflectbox{$\fix@cev{#1}{-}\m@th#1#2\fix@cev{#1}{+}$}}$}%
	\fix@cev{#1}{-}%
}
\newcommand{\fix@cev}[2]{%
	\ifx#1\displaystyle
	\mkern#23mu
	\else
	\ifx#1\textstyle
	\mkern#23mu
	\else
	\ifx#1\scriptstyle
	\mkern#22mu
	\else
	\mkern#22mu
	\fi
	\fi
	\fi
}
\usepackage{soul}
\usepackage{mathtools}  % Agregar en el preámbulo
\usepackage{pgfplots}
\usepackage{tikz}
\usepackage{ifthen}
\usetikzlibrary{shapes.geometric, positioning, calc}
\tikzset{%
    mynode/.style={%
        circle, minimum size=1.21mm, inner sep=0pt, draw=black, fill=black
    }
}

\usepackage{todonotes}

\definecolor{lightblue}{rgb}{0.68,0.85,0.9}

\begin{document}
\title{Isometric and induced path partitions: a new upper bound and a characterization of some extremal graphs}

\author{Irena Penev\thanks{Computer Science Institute (I\'UUK, MFF), Charles University, Prague, Czech Republic. Supported by GA\v{C}R grant 25-17377S and by MSCA-RISE-2020-101007705 project \textit{RandNET}. Email: {ipenev@iuuk.mff.cuni.cz}.  } ~~ R. B. Sandeep\thanks{Department of Computer Science and Engineering, Indian Institute of Technology Dharwad, India. Supported by ANRF (SERB) grant MTR/2022/000692. Email: {sandeeprb@iitdh.ac.in. }} ~~ D.K. Supraja\thanks{Indian Institute of Technology Dharwad, India. Email: {dksupraja95@gmail.com}. } ~~
S Taruni\thanks{Centro de Modelamiento Matemático (CNRS IRL2807), Universidad de Chile, Santiago, Chile. Supported by Centro de Modelamiento Matemático (CMM) BASAL fund FB210005 for center of excellence from ANID-Chile and by MSCA-RISE-2020-101007705 project \textit{RandNET}. Email: {tsridhar@cmm.uchile.cl}. } ~~ }
\date{}
\maketitle
\begin{abstract}
\noindent 
    An \textit{isometric path} is a shortest path between two vertices. An \textit{isometric path partition} (IPP) of a graph $G$ is a set $\mathcal{I}$ of vertex-disjoint isometric paths in $G$ that partition the vertices of $G$. The \textit{isometric path partition number} of $G$, denoted by $\text{ipp}(G)$, is the minimum cardinality of an IPP of~$G$. An \textit{induced path partition} (IndPP) of a graph $G$ is a set $\mathcal{I}$ of vertex-disjoint induced paths in~$G$ that partition the vertices of $G$. The \textit{induced path partition number} of $G$, denoted by $\text{indpp}(G)$, is the minimum cardinality of an IndPP of $G$. In this article, we study both these parameters and observe that every graph $G$ satisfies $\text{indpp}(G) \leq \text{ipp}(G) \leq |V(G)| - \nu(G)$, where $\nu(G)$ is the matching number of $G$. We further prove that a connected graph $G$ is extremal with respect to this upper bound, i.e.\ satisfies $\text{ipp}(G) = |V(G)| - \nu(G)$, (resp.\ $\text{indpp}(G) = |V(G)| - \nu(G)$), if and only if either (i) all blocks of $G$ are odd complete graphs, or (ii) all blocks of $G$ except one are odd complete graphs, and the unique block $B$ of $G$ that is not an odd complete graph is even and satisfies $\text{ipp}(B) = |V(B)| - \nu(B)$ (resp.\ $\text{indpp}(B) = |V(B)| - \nu(B)$). As corollaries of these results, we obtain a full structural characterization of all connected odd graphs that are extremal with respect to our upper bound, as well as of all extremal block graphs. 
%We further characterize all connected graphs for which this upper bound is attained. As an application of this characterization, we identify all connected graphs of odd order, even order, and all block graphs that achieve this bound.

\end{abstract}

% \usepackage{fullpage}
% %\usepackage{geometry}
% %\geometry{
% 	%	left=40mm,
% 	%	right=40mm,
% 	%	top=35mm,
% 	%	bottom=35mm, 
% 	%}
% \usepackage{xcolor}
% \usepackage{amssymb}
% \usepackage{amsmath, amsthm, amsfonts}
% \usepackage{float}
% \usepackage{graphicx}
% \usepackage{lineno}
% \usepackage{enumerate}
% \usepackage{longtable}
% \usepackage{environ}
% \usepackage{bm}
% \usepackage{thm-restate}
% \usepackage{hyperref}
% \usepackage{cleveref,lipsum}

% \usepackage[numbers]{natbib}
% \usepackage{xurl}

% \usepackage{soul}

\

\begin{center} 
\fbox{\begin{minipage}{35em}
\noindent
This manuscript is a slightly improved version of the following paper: 
\begin{quote} 
Irena Penev,  R.B.~Sandeep, D.K.~Supraja, S~Taruni, ``Isometric and induced path partitions: a new upper bound and a characterization of some extremal graphs''. {\em Discrete Applied Mathematics}, 392 (2026), 260--273. 
\end{quote} 
Sections~1 and~4 (i.e.\ the introduction and the concluding remarks) differ from the published version. The most significant difference is that Definition~1.3 and Proposition~1.6 from this manuscript are missing from the published version, and moreover, Proposition~1.4 from this manuscript (which corresponds to Proposition~1.3 of the published paper) also includes $\text{pp}$ (rather than just $\text{indpp}$ and $\text{ipp}$). The main theorems of this manuscript are Theorems~1.8 and~1.9, and they are exactly the same as Theorems~1.6 and~1.7, respectively, of the published version. Moreover, we have corrected a few typos and minor errors that appeared in section~2; a full list of these corrections can be found at the end of the manuscript (after the references). 
\end{minipage}}
\end{center} 

\bigskip 
\bigskip 

\section{Introduction} 
In what follows, all graphs are finite, simple, and nonnull. The vertex set and the edge set of a graph $G$ are denoted by $V(G)$ and $E(G)$, respectively, and we say that the graph $G$ is {\em even} (resp.\ {\em odd}) if $|V(G)|$ is even (resp.\ odd). An {\em $n$-vertex graph} is a graph that contains exactly $n$ vertices. 

Graph partitioning and covering problems are central topics in graph theory and algorithms, including problems such as dominating set (covering by stars), clique covering (covering by cliques), coloring (partitioning by independent sets), and covering/partitioning by paths or cycles. In particular, the problems of partitioning and covering by paths have their connections to important graph-theoretic results such as the Gallai-Milgram theorem~\cite{gallai1960verallgemeinerung} and Berge's path partition conjecture~\cite{berge1983path}, and also to the celebrated Graph Minor project~\cite{kawarabayashi2012disjoint,robertson1995graph}. These problems also have applications in various domains, including artificial intelligence~\cite{guo2018fast}, bioinformatics~\cite{lafond2025path}, code optimization~\cite{boesch1977covering}, program testing~\cite{ntafos1979path}, to name a few. In recent years, various types of path partitions have been studied, such as unrestricted path partition~\cite{corniel2013,garey1976,skupien1974path}, induced path partition~\cite{chartrand1994, LE2003, PAN2007121}, and isometric (i.e.\ shortest) path partition~\cite{chakraborty2026covering,fernau2025parameterizing}. See~\cite{fernau2025parameterizing, manuel2018revisiting} for a broader overview of path partition and related problems in the literature. %We note that all isometric paths are induced, but not all induced paths are isometric. 

The focus of this article is on the problem of partitioning the vertex set of a graph into isometric paths. Let us be more precise. A {\em path} (resp.\ {\em induced path}) in a graph $G$ is a subgraph (resp.\ induced subgraph) $P$ of $G$ such that $P$ is a path. The {\em length} of a path $P$, denoted by $\text{length}(P)$, is the number of edges that it contains. Given vertices $x$ and $y$ in a graph $G$, a {\em shortest} or {\em isometric} path between $x$ and~$y$ in $G$ is a path in $G$ whose endpoints are $x$ and $y$, and is of smallest possible length among all paths between $x$ and $y$ in $G$. When we say that ``$P$ is an isometric path in $G$,'' we mean that $P$ is a shortest path between its endpoints in $G$ (here, it is possible that $P$ has only one vertex).  We now state the main definitions.

\begin{definition}
    A {\em path partition} (or {\em PP} for short) of a graph $G$ is a collection of paths in $G$ whose vertex sets partition $V(G)$. The {\em path partition number} of $G$, denoted by $\text{pp}(G)$, is the smallest cardinality of any PP of $G$. A {\em minimum} PP of $G$ is a PP of $G$ of cardinality $\text{pp}(G)$. 
\end{definition}

\begin{definition}
    An {\em induced path partition} (or {\em IndPP} for short) of a graph $G$ is a collection of induced paths in $G$ whose vertex sets partition $V(G)$. The {\em induced path partition number} of $G$, denoted by $\text{indpp}(G)$, is the smallest cardinality of any IndPP of $G$. A {\em minimum} IndPP of $G$ is an IndPP of $G$ of cardinality $\text{indpp}(G)$. 
\end{definition}

\begin{definition}
    An {\em isometric path partition} (or {\em IPP} for short) of a graph $G$ is a collection of isometric paths in $G$ whose vertex sets partition $V(G)$. The {\em isometric path partition number} of $G$, denoted by $\text{ipp}(G)$, is the smallest cardinality of any IPP of $G$. A {\em minimum} IPP of $G$ is an IPP of $G$ of cardinality $\text{ipp}(G)$. 
\end{definition}

Clearly, every isometric path in~$G$ is an induced path of~$G$, but the converse is not true in general. We note, however, that any induced path of length at most two in~$G$ is an isometric path of~$G$, a fact that we will repeatedly use in this paper. We further observe that any IPP of a graph $G$ is an IndPP of~$G$, and that any IndPP of~$G$ is a PP of~$G$; consequently, 
\begin{displaymath} 
\begin{array}{ccc cc} 
\text{pp}(G) & \leq & \text{indpp}(G) & \leq & \text{ipp}(G). 
\end{array} 
\end{displaymath}

The {\em distance} between vertices $x$ and $y$ in a graph $G$, denoted by $\text{dist}_G(x,y)$, is defined as the length of a shortest path between $x$ and $y$ (if no such path exists, i.e.\ if $x$ and $y$ belong to different connected components of $G$, then we define $\text{dist}_G(x,y) := \infty$); if $x = y$, then clearly, $\text{dist}_G(x,y) = 0$. The {\em diameter} of a graph $G$ is defined as $\text{diam}(G) := \text{max}\big\{\text{dist}_G(x,y) \mid x,y \in V(G)\big\}$. Clearly, any isometric path in~$G$ is of length at most $\text{diam}(G)$ (i.e.\ contains at most $\text{diam}(G)+1$ many vertices), and consequently, as pointed out in~\cite{PManuel2021}, every connected graph $G$ satisfies
\begin{displaymath} 
\begin{array}{rcl} 
\text{ipp}(G) & \geq & \Big\lceil \frac{|V(G)|}{\text{diam}(G)+1} \Big\rceil. 
\end{array} 
\end{displaymath}

In this article, we give an upper bound for the PP, IndPP, and IPP numbers in terms of the matching number; ``PP-extremal'' (resp.\ ``IndPP-extremal,'' ``IPP-extremal) graphs are those whose PP number (resp.\ IndPP number, IPP number) equals this upper bound. PP-extremal graphs are easy to characterize (see Proposition~\ref{prop-PP-extremal}). The main result of this article is a characterization of the connected IPP-extremal and IndPP-extremal graphs in terms of their blocks. These results are described in subsection~\ref{subsec:contribution} (see Proposition~\ref{prop-upper-bound-trivial}, Theorem~\ref{thm-main} and Theorem~\ref{thm-main-indpp}). As a corollary of our main results (Theorems~\ref{thm-main} and~\ref{thm-main-indpp}), we give a full structural characterization of the connected odd IndPP-extremal and IPP-extremal graphs (see Corollary~\ref{cor-main-odd}), as well of the IndPP-extremal and IPP-extremal block graphs (see Corollary~\ref{cor-extremal-block-graph}). 

The remainder of the paper is organized as follows. In subsection~\ref{subsec:previous}, we provide the context for our work. In subsection~\ref{subsec:contribution}, we describe our main results. In subsection~\ref{subsec:notation}, we introduce some (mostly standard) terminology and notation that we will use throughout the paper. In section~\ref{sec:proof-main}, we prove Theorems~\ref{thm-main} and~\ref{thm-main-indpp} (our main theorems). In section~\ref{sec:corollaries}, we derive some corollaries of Theorems~\ref{thm-main} and~\ref{thm-main-indpp}. Finally, in section~\ref{sec:conclusion}, we propose some open questions.

\subsection{Previous work} \label{subsec:previous}

For a set $X$, a collection $\mathcal{Z}$ of subsets of $X$ is said to {\em cover} $X$ if for all $x \in X$, there exists some $Z \in \mathcal{Z}$ such that $x \in Z$. An {\em isometric path cover} (or {\em IPC} for short) of a graph $G$ is a collection of isometric paths of $G$ whose vertex sets cover $V(G)$. The {\em isometric path cover number} of $G$, denoted by $\text{ipc}(G)$, is the minimum cardinality of any IPC of $G$. IPC was introduced by Fitzpatrick~\cite{fitzpatrick1998aspects}, and it is obviously related to IPP. Indeed, any IPP of a graph $G$ is, in particular, an IPC of $G$, and it follows that $\text{ipp}(G) \geq \text{ipc}(G)$. 

Fitzpatrick referred to the IPC number by different names such as the \textit{precinct number}~\cite{fitzpatrick1998aspects} or the \textit{isometric path number}~\cite{FITZPATRICK2001253}. This concept proved useful in a particular variation of the game of cops and robbers. For any graph $G$ where a single cop suffices to capture the robber, a strategy for the cop can be determined by projecting the game positions onto isometric paths within the graph. Fitzpatrick~\cite{fitzpatrick1998aspects} observed that all graphs $G$ satisfy $\text{ipc}(G) \geq \big\lceil |V(G)|/(\text{diam}(G)+1) \big\rceil$, and he described some classes of graphs that reach this lower bound (i.e.\ graphs for which this inequality becomes an equality). In particular, Fitzpatrick et al.~\cite{FITZPATRICK2001253} proved that the hypercube $Q_n$ satisfies $\text{ipc}(Q_n) \geq 2^n / (n+1)$, and that this lower bound is in fact tight when $n+1$ is a power of $2$. 

The exact values of the IPP number have been established for hypercubes $Q_n$ when $n+1$ is a power of~$2$. For $r \times r$ tori, the value is $r$ when $r$ is even and is known to lie in $\{r, r+1\}$ when $r$ is odd~\cite{PManuel2021}. We remark that the results of~\cite{PManuel2021} relied on previously established results for the IPC number, as well as on the above-mentioned fact that every graph $G$ satisfies $\text{ipp}(G) \geq \text{ipc}(G)$. To our knowledge, no other significant structural results have been obtained for either IPP or IndPP. However, there have been a number of algorithmic results for these two parameters. In particular, Le et al.~\cite{LE2003} proved that it is NP-complete to decide whether or not the vertex set of a connected graph can be partitioned into two subsets each of which induces a path; consequently, computing the IndPP number of a graph is NP-hard. Similarly, Manuel~\cite{PManuel2021} showed that computing an optimal IPP is NP-hard on undirected graphs; this was extended in~\cite{fernau2025parameterizing} to prove that IPP (studied under the name \textit{shortest path partition}, SPP in \cite{fernau2025parameterizing}) is NP-hard even when restricted to bipartite graphs that are sparse (in particular, when they have degeneracy at most 5). It was further shown in~\cite{fernau2025parameterizing} that IndPP and IPP belong to FPT when parameterized by standard structural parameters such as the vertex cover and the neighborhood diversity (this result was obtained using Integer Linear Programming), and moreover, that IPP can be solved in XP-time on undirected graphs when parameterized by solution size. IPP is known to be NP-hard on split graphs while the problem is polynomial-time solvable for cographs and chain graphs~\cite{chakraborty2026covering}. Additionally, polynomial-time algorithms are known for computing optimal (unrestricted) path partitions in trees~\cite{skupien1974path} and for computing optimal induced path partitions in graphs with special blocks~\cite{PAN2007121}. Since, for trees, the (unrestricted) path partition and induced path partition coincide with the isometric path partition, this immediately yields polynomial-time algorithms for computing a minimum IPP and IndPP in trees~\cite{skupien1974path}.

\subsection{Our contribution} \label{subsec:contribution} 

A {\em matching} of a graph $G$ is a set of edges of $G$, no two of which share an endpoint. The {\em matching number} of $G$, denoted by $\nu(G)$, is the largest cardinality of any matching in $G$; clearly, $|V(G)| \geq 2\nu(G)$. A {\em maximum matching} of $G$ is a matching of $G$ of cardinality $\nu(G)$. For a vertex $v$ and a matching $M$ of $G$, we say that $v$ is {\em $M$-saturated} if $v$ is incident with some edge in $M$, and otherwise, we say that $v$ is {\em $M$-unsaturated}. A {\em perfect matching} of $G$ is a matching $M$ of $G$ such that every vertex of $G$ is $M$-saturated; clearly, $G$ has a perfect matching if and only if $|V(G)| = 2\nu(G)$. 

For notational convenience, we will identify any one-vertex path with its unique vertex, and we will identify any two-vertex path with its unique edge. Clearly, any one- or two-vertex path is isometric. So, if $M$ is a matching of a graph $G$, and $U$ is the set of all $M$-unsaturated vertices of $G$, then $\mathcal{I} := M \cup U$ is an IPP of $G$, and we see that $\text{ipp}(G) \leq |\mathcal{I}| = |M|+|U| = |M|+\big(|V(G)|-2|M|\big) = |V(G)|-|M|$. By choosing $M$ to be a maximum matching of $G$, so that $|M| = \nu(G)$, we obtain the inequality $\text{ipp}(G) \leq |V(G)|-\nu(G)$, an upper bound on the IPP number. Combined with our previously mentioned inequality between the PP, IndPP, and IPP numbers, this yields the following proposition. 

\begin{proposition} \label{prop-upper-bound-trivial} Every graph $G$ satisfies $\text{pp}(G) \leq \text{indpp}(G) \leq \text{ipp}(G) \leq |V(G)|-\nu(G)$. 
\end{proposition}

Let us say that a graph $G$ is {\em IPP-extremal} if it satisfies $\text{ipp}(G) = |V(G)|-\nu(G)$ (so that the rightmost inequality from Proposition~\ref{prop-upper-bound-trivial} becomes an equality). Similarly, a graph $G$ is {\em IndPP-extremal} if it satisfies $\text{indpp}(G) = |V(G)|-\nu(G)$, and $G$ is {\em PP-extremal} if it satisfies $\text{pp}(G) = |V(G)|-\nu(G)$. By Proposition~\ref{prop-upper-bound-trivial}, every PP-extremal graph is IndPP-extremal, and every IndPP-extremal graph is IPP-extremal; the latter fact is also stated below (see Proposition~\ref{prop-indpp-implies-ipp}) for future reference. There are, however, IndPP-extremal graphs that are not PP-extremal; examples include complete graphs on more than two vertices, as well as the cycle $C_4$. Likewise, not all IPP-extremal graphs are IndPP-extremal; for example, it is easy to see that the wheel on six vertices represented in Figure~\ref{fig:indpp} is IPP-extremal, but not IndPP-extremal. 

\begin{figure}[h]
    \centering
   \begin{tikzpicture}
\centering
     \node[minimum size=2.2cm, regular polygon, regular polygon sides=5, rotate=0] (pent) at (11.25,-1.7) {};
        \foreach \x in {1,2,...,5}{%
            \node[mynode] at (pent.corner \x) (f\x) {};
    }
    \draw (f1)to(f2)to(f3)to(f4)to(f5)to(f1);
    \node[mynode] (a) at (11.25,-1.7) {};
        \foreach \x in {1,2,...,5}{%
            \draw (f\x)to(a);
            }
\end{tikzpicture}
    \caption{A graph that is IPP-extremal, but not IndPP-extremal.}
    \label{fig:indpp}
\end{figure}

\begin{proposition} \label{prop-indpp-implies-ipp} If a graph $G$ is IndPP-extremal, then both the following hold: 
\begin{itemize} 
\item $G$ is IPP-extremal; 
\item every minimum IPP of $G$ is also a minimum IndPP of $G$. 
\end{itemize} 
\end{proposition} 
\begin{proof} 
Fix an IndPP-extremal graph $G$, so that $\text{indpp}(G) = |V(G)|-\nu(G)$. Combined with Proposition~\ref{prop-upper-bound-trivial}, this implies that $|V(G)|-\nu(G) = \text{indpp}(G) \leq \text{ipp}(G) \leq |V(G)|-\nu(G)$, and consequently, $\text{indpp}(G) = \text{ipp}(G) = |V(G)|-\nu(G)$. Thus, $G$ is IPP-extremal. Since every IPP of $G$ is an IndPP of $G$, it further follows that any minimum IPP of $G$ is a minimum IndPP of $G$. 
\end{proof}

PP-extremal graphs are simple to characterize. Indeed, as pointed out to us by an anonymous referee, PP-extremal graphs are precisely those that are disjoint unions of $K_1$'s and $K_2$'s, i.e.\ those that have no components on more than two vertices. More precisely, we have the following proposition. 

\begin{proposition} \label{prop-PP-extremal} For every graph $G$, the following are equivalent: 
\begin{itemize} 
\item $G$ is PP-extremal; 
\item $\text{pp}(G) = \text{indpp}(G) = \text{ipp}(G) = |V(G)|-\nu(G)$; 
\item $G$ has no components on more than two vertices. 
\end{itemize} 
\end{proposition} 
\begin{proof} 
Fix a graph $G$. The equivalence of the first two statements follows from Proposition~\ref{prop-upper-bound-trivial} and from the definition of a PP-extremal graph. Moreover, it is clear that the third statement implies the second. It now remains to show that the first statement implies the third. We prove the contrapositive: we assume that $G$ has a component, call it $C$, on at least three vertices, and we show that $G$ has a PP of cardinality $|V(G)|-\nu(G)-1$, so that $G$ is not PP-extremal. 

Let $M$ be a maximum matching of $G$, and let $U$ be the set of all $M$-unsaturated vertices of $G$, so that $\mathcal{P} := M \cup U$ is a PP of $G$ of cardinality $|V(G)|-\nu(G)$. Clearly, the component $C$ contains at least one edge of $M$ (otherwise, we consider any edge $e$ of $C$, and we observe that $M \cup \{e\}$ is a matching of $G$, contrary to the maximality of $M$); let $uv \in E(C) \cap M$. Since $C$ is connected and has at least three vertices, we may assume by symmetry that there is a vertex $w \in V(C) \setminus \{u,v\}$ such that $vw \in E(C)$. If $w$ is $M$-unsaturated, then $\mathcal{P}' := \big(M \setminus \{uv\}\big) \cup \big(U \setminus \{w\}\big) \cup \{uvw\}$ is a PP of $G$ of cardinality $|\mathcal{P}|-1 = |V(G)|-\nu(G)-1$, and we are done. We may therefore assume that $w$ is $M$-saturated, so that there exists a vertex $w' \in V(C) \setminus \{u,v,w\}$ such that $ww' \in M$. But now $\mathcal{P}'' := \big(M \setminus \{uv,ww'\}\big) \cup U \cup \{uvww'\}$ is a PP of $G$ of cardinality $|\mathcal{P}|-1 = |V(G)|-\nu(G)-1$, and again we are done. 
\end{proof}

%Since every graph $G$ satisfied $\text{indpp}(G) \leq \text{ipp}(G)$, we see that every IndPP-extemal graph is IPP-extremal. We note, however, that the converse does not hold. For example, the wheel on six vertices (see Figure~\ref{fig:indpp}) is IPP-extremal, but is not IndPP-extremal. 

In the remainder of the paper, we turn our attention to IndPP-extremal and IPP-extremal graphs. Our goal is to give a partial characterization of IPP-extremal graphs and IndPP-extremal graphs. The following proposition reduces this problem to the connected case. 

\begin{proposition} \label{prop-components} 
Let $G$ be a graph. Then both the following hold: 
\begin{enumerate}[(a)] 
\item $G$ is IPP-extremal if and only if all components of $G$ are IPP-extremal; 
\item $G$ is IndPP-extremal if and only if all components of $G$ are IndPP-extremal. 
\end{enumerate} 
\end{proposition}
\begin{proof} 
We prove~(a); the proof of~(b) is completely analogous. Let $G_1,\dots,G_t$ be the connected components of $G$. For each $i \in \{1,\dots,t\}$, let $M_i$ be a maximum matching of $G_i$, and let $\mathcal{I}_i$ be a minimum IPP of $G_i$. Then $M := M_1 \cup \dots \cup M_t$ is a maximum matching of $G$, and $\mathcal{I} := \mathcal{I}_1 \cup \dots \cup \mathcal{I}_t$ is a minimum IPP of $G$. We now have that 
\begin{longtable}{rcl} 
$\text{ipp}(G)$ & $=$ & $\sum\limits_{i=1}^t \text{ipp}(G_i)$ 
\\
\\
& $\stackrel{(*)}{\leq}$ & $\sum\limits_{i=1}^t \Big(|V(G_i)|-\nu(G_i)\Big)$ 
\\
\\
& $=$ & $\Big(\sum\limits_{i=1}^t |V(G_i)|\Big)-\Big(\sum\limits_{i=1}^t \nu(G_i)\Big)$ 
\\
\\
& $=$ & $|V(G)|-\nu(G),$ 
\end{longtable} 
\noindent 
where (*) follows from Proposition~\ref{prop-upper-bound-trivial}. Moreover, the inequality (*) is an equality precisely when $\text{ipp}(G_i) = |V(G_i)|-\nu(G_i)$ for all $i \in \{1,\dots,t\}$. This proves~(a). 
\end{proof}

A {\em block} of a graph is a maximal biconnected induced subgraph of that graph. Our main results are the following.

\begin{theorem} \label{thm-main} Let $G$ be a connected graph. Then $G$ is IPP-extremal if and only if $G$ satisfies one of the following: 
\begin{enumerate}[(i)] 
\item all blocks of $G$ are odd complete graphs; 
\item all blocks of $G$ except one are odd complete graphs, and the unique block of $G$ that is not an odd complete graph is even and IPP-extremal. 
\end{enumerate} 
\end{theorem} 

\begin{theorem} \label{thm-main-indpp} Let $G$ be a connected graph. Then $G$ is IndPP-extremal if and only if $G$ satisfies one of the following: 
\begin{enumerate}[(i)] 
\item all blocks of $G$ are odd complete graphs; 
\item all blocks of $G$ except one are odd complete graphs, and the unique block of $G$ that is not an odd complete graph is even and IndPP-extremal. 
\end{enumerate} 
\end{theorem}

\subsection{Terminology and notation} \label{subsec:notation} 

When we write ``$P = p_1\dots p_t$ is a path'' ($t \geq 1$) we mean that $P$ is a $t$-vertex graph with $V(P) = \{p_1,\dots,p_t\}$ and $E(P) = \{p_ip_{i+1} \mid 1 \leq i \leq t-1\}$; under these circumstances, we say that $p_1$ and $p_t$ are the {\em endpoints} of the path $P$ (or that $P$ is a path {\em between} $p_1$ and $p_t$), whereas $p_2,\dots,p_{t-1}$ are the {\em internal vertices} of $P$ (if $t = 1$, then $p_1$ is the only endpoint of $P$, and $P$ has no internal vertices). 

For a vertex $v$ in a graph $G$, the {\em open neighborhood} of $v$ in $G$, denoted by $N_G(v)$, is the set of all neighbors of $v$ in $G$, whereas the {\em closed neighborhood} of $v$ in $G$ is the set $N_G[v] = \{v\} \cup N_G(v)$. For a vertex $v \in V(G)$ and a set $S \subseteq V(G) \setminus \{v\}$, we say that $v$ is {\em complete} (resp.\ {\em anticomplete}) to $S$ in $G$ if $v$ is adjacent (resp.\ nonadjacent) to all vertices of $S$, and we say that $v$ is {\em mixed} on $S$ in $G$ if $v$ is neither complete nor anticomplete to $S$ in $G$ (i.e.\ if $v$ has both a neighbor and a nonneighbor in $S$). Recall that, formally, an edge is simply a set of two distinct vertices (the endpoints of the edge). So, it makes sense to speak of a vertex being (anti)complete to or mixed on an edge (as long as the vertex is not incident with the edge). For disjoint sets $X,Y \subseteq V(G)$, we say that $X$ is {\em complete} (resp.\ {\em anticomplete}) to $Y$ in~$G$ if every vertex in $X$ is complete (resp.\ anticomplete) to $Y$ in $G$. 

A {\em clique} (resp.\ {\em stable set}) in a graph $G$ is any set of pairwise adjacent (resp.\ nonadjacent) vertices of~$G$. Note that $\emptyset$ is both a clique and a stable set in any graph. A clique or a stable set is {\em even} (resp.\ {\em odd}) if its cardinality is even (resp.\ odd). A vertex $v$ of a graph $G$ is {\em simplicial} if $N_G(v)$ is a (possibly empty) clique of $G$. 

For a graph $G$ and a nonempty set $S \subseteq V(G)$, we denote by $G[S]$ the subgraph of $G$ induced by $S$. For $S \subsetneqq V(G)$, we denote by $G \setminus S$ the subgraph of $G$ obtained by deleting all vertices of $S$, that is, $G \setminus S := G[V(G) \setminus S]$; if $S = \{v\}$, we may write $G \setminus v$ instead of $G \setminus S$ or $G \setminus \{v\}$. For a graph $G$, a {\em cutset} of $G$ is any set $S \subsetneqq V(G)$ such that $G \setminus S$ contains more connected components than $G$ does. For a graph $G$ on at least two vertices, a {\em cut-vertex} of $G$ is a vertex $v$ of $G$ such that $G \setminus v$ has more connected components than $G$ has (i.e.\ such that $\{v\}$ is a cutset of $G$). 

A graph $G$ is {\em biconnected} if it is connected and contains no cut-vertices. Note that all complete graphs are biconnected. A {\em block} of a graph $G$ is a maximal biconnected induced subgraph of $G$. Consistently with our terminology above, a block $B$ of a graph $G$ is {\em even} (resp.\ {\em odd}) if $|V(B)|$ is even (resp.\ odd). Clearly, every connected graph $G$ that is not biconnected contains a cut-vertex. A {\em block graph} is a graph, all of whose blocks are complete graphs. A {\em leaf-block} of a graph $G$ is a block of $G$ that contains exactly one cut-vertex of $G$. Note that any connected graph that is not biconnected contains at least two leaf-blocks. A {\em leaf-clique} of a graph $G$ is a nonempty clique $C$ of $G$ such that $G[C]$ is a leaf-block of~$G$.

\section{Proof of Theorems~\ref{thm-main} and~\ref{thm-main-indpp}} \label{sec:proof-main}

\begin{proposition} \label{prop-no-unsaturated-mixed-on-matching-edge} Let $G$ be a connected IPP-extremal graph, and let $M$ be a maximum matching of $G$. Then no $M$-unsaturated vertex of $G$ is mixed on any edge of $M$. Moreover, $G$ contains at most one $M$-unsaturated vertex. 
\end{proposition} 
\begin{proof} 
Let $U$ be the set of all $M$-unsaturated vertices of $G$. We must show that no vertex in $U$ is mixed on any edge in $M$, and that $|U| \leq 1$. Since $M$ is a maximum matching of $G$, we have that $|M| = \nu(G)$, and consequently, $|U| = |V(G)|-2|M| = |V(G)|-2\nu(G)$. So, $\mathcal{I} := M \cup U$ is an IPP of $G$ of size $|V(G)|-\nu(G)$; since $\text{ipp}(G) = |V(G)|-\nu(G)$, it follows that $\mathcal{I}$ is in fact a minimum IPP of $G$. 

First, suppose toward a contradiction that some $M$-unsaturated vertex $u$ is mixed on an edge $xy \in M$. By symmetry, we may assume that $ux \in E(G)$ and $uy \notin E(G)$. But now, $\big(M \setminus \{xy\}\big) \cup \big(U \setminus \{u\}\big) \cup \{uxy\}$ is an IPP of $G$ of size $\big(|M|-1\big)+\big(|U|-1\big)+1 = |V(G)|-\nu(G)-1$, contrary to the assumption that $G$ is IPP-extremal and therefore satisfies $\text{ipp}(G) = |V(G)|-\nu(G)$. 

It remains to show that $|U| \leq 1$. Suppose otherwise, and fix distinct $u_1,u_2 \in U$. Let $P'$ be a shortest path between $u_1$ and $u_2$ in $G$. But then 
\begin{displaymath} 
\begin{array}{rcl} 
\mathcal{I'} & := & \{P'\} \cup \big\{P \in \mathcal{I} \mid V(P') \cap V(P) = \emptyset \big\} \cup \big\{ v \in V(G) \setminus V(P') \mid \text{$\exists v' \in V(P')$ s.t.\ $vv' \in M$}\big\}
\end{array} 
\end{displaymath} 
is an IPP of $G$ size at most $|\mathcal{I}|-1$, contrary to the minimality of $\mathcal{I}$. 
\end{proof} 

\begin{proposition} \label{prop-basic-even} Let $G$ be a connected even IPP-extremal graph. Then $|V(G)| = 2\nu(G)$ and $\text{ipp}(G) = \nu(G)$. Moreover, $G$ admits a perfect matching, and any perfect matching of $G$ is a minimum IPP of $G$. 
\end{proposition} 
\begin{proof} 
Fix a maximum matching $M$ of $G$. By Proposition~\ref{prop-no-unsaturated-mixed-on-matching-edge}, at most one vertex of $G$ is $M$-unsaturated. But since $G$ is even, this implies that no vertex of $G$ is $M$-unsaturated, that is, $M$ is a perfect matching of $G$. Consequently, $|V(G)| = 2\nu(G)$. Since $\text{ipp}(G) = |V(G)|-\nu(G)$, it follows that $\text{ipp}(G) = \nu(G)$. Finally, since any perfect matching is an IPP of $G$ of size $\nu(G)$, it follows that any perfect matching of $G$ is a minimum IPP of $G$. 
\end{proof}

\begin{proposition} \label{prop-basic-odd} Let $G$ be a connected odd IPP-extremal graph. Then $|V(G)| = 2\nu(G)+1$ and $\text{ipp}(G) = \nu(G)+1$.  Moreover, for every vertex $u \in V(G)$, the graph $G \setminus u$ admits a perfect matching, and furthermore, for any perfect matching $M_u$ of $G \setminus u$, we have that $M_u \cup \{u\}$ is a minimum IPP of $G$. 
\end{proposition} 
\begin{proof} 
Since $G$ is odd, Proposition~\ref{prop-no-unsaturated-mixed-on-matching-edge} guarantees that for every maximum matching $M$ of $G$, exactly one vertex of $G$ is $M$-unsaturated; consequently, $|V(G)| = 2\nu(G)+1$. Since $\text{ipp}(G) = |V(G)|-\nu(G)$, it follows that $\text{ipp}(G) = \nu(G)+1$. 

It is clear that if $M_u$ is a perfect matching of $G \setminus u$ for some vertex $u$, then $M_u$ is a maximum matching of $G$, and that $M_u \cup \{u\}$ is an IPP of $G$; since $\text{ipp}(G) = \nu(G)+1$, we see that $M_u \cup \{u\}$ is in fact a minimum IPP of $G$. 

It remains to prove that for all $u \in V(G)$, the graph $G \setminus u$ admits a perfect matching. Set 
\begin{displaymath} 
\begin{array}{rcl} 
\mathcal{U} & := & \{u \in V(G) \mid \text{$G \setminus u$ admits a perfect matching}\}. 
\end{array} 
\end{displaymath} 
We must show that $\mathcal{U} = V(G)$. Suppose otherwise, so that $V(G) \setminus \mathcal{U} \neq \emptyset$. Since $|V(G)| = 2\nu(G)+1$, it is clear that $\mathcal{U} \neq \emptyset$.\footnote{Indeed, if $M$ is any maximum matching of $G$, then the fact that $|V(G)| = 2\nu(G)+1$ implies that exactly one vertex of $G$, call it $v_M$, is $M$-unsaturated; but clearly, $M$ is a perfect matching of $G \setminus v_M$, and so $v_M \in \mathcal{U}$.} Since $G$ is connected, there exist adjacent vertices $u \in \mathcal{U}$ and $v \in V(G) \setminus \mathcal{U}$. Since $u \in \mathcal{U}$, $G \setminus u$ admits a perfect matching, call it $M_u$. Clearly, $M_u$ is a maximum matching of $G$, and $u$ is the only $M_u$-unsaturated vertex of $G$. Then $v$ is $M_u$-saturated, and therefore, there exists some $v' \in V(G) \setminus \{u,v\}$ such that $vv' \in M_u$. Since $u$ is $M_u$-unsaturated and adjacent to $v$, Proposition~\ref{prop-no-unsaturated-mixed-on-matching-edge} guarantees that $u$ is also adjacent to $v'$.\footnote{Otherwise, the $M_u$-unsaturated vertex $u$ would be mixed on the matching edge $vv' \in M_u$, contrary to Proposition~\ref{prop-no-unsaturated-mixed-on-matching-edge}.} But now $M_v := \big(M_u \setminus \{vv'\}\big) \cup \{uv'\}$ is a perfect matching of $G \setminus v$, contrary to the fact that $v \notin \mathcal{U}$. 
\end{proof}

\begin{proposition} \label{prop-delete-edge} Let $C$ be a leaf-clique of a graph $G$, let $v$ be the (unique) cut-vertex of $G$ that belongs to $C$, and let $x$ and $y$ be distinct vertices in $C \setminus \{v\}$. Set $G' := G \setminus \{x,y\}$. Then both the following hold: 
\begin{enumerate}[(a)] 
\item $G$ is IPP-extremal if and only if $G'$ is IPP-extremal. 
\item $G$ is IndPP-extremal if and only if $G'$ is IndPP-extremal. 
\end{enumerate} 
\end{proposition} 
\begin{proof} 
We prove (a) and (b) simultaneously. We begin by observing all vertices $z \in C \setminus \{v\}$ satisfy $N_G[z] = C$, and consequently, all vertices in $C \setminus \{v\}$ are simplicial in $G$. In particular, we have that $N_G[x] = N_G[y] = C$, and that vertices $x$ and $y$ are simplicial in $G$. 

\begin{adjustwidth}{1cm}{1cm} 
\begin{claim} \label{prop-delete-edge-claim-delete-edge} All the following hold: 
\begin{enumerate}[(1)] 
\item $|V(G')| = |V(G)|-2$; 
\item $\nu(G') = \nu(G)-1$; 
\item $|V(G')|-\nu(G') = |V(G)|-\nu(G)-1$. 
\end{enumerate} 
\end{claim} 
\end{adjustwidth} 
\noindent 
{\em Proof of Claim~\ref{prop-delete-edge-claim-delete-edge}.} 
Part~(1) follows immediately from the construction, and part~(3) follows immediately from~(1) and~(2). Thus, we just need to prove~(2). 

Clearly, if $M'$ is a maximum matching of $G'$, then $M' \cup \{xy\}$ is a matching of $G$, and so $\nu(G) \geq |M' \cup \{xy\}| = |M'|+1 = \nu(G')+1$, which implies that $\nu(G') \leq \nu(G)-1$. 

It remains to show that $\nu(G') \geq \nu(G)-1$. Clearly, it is enough to show that $G$ contains a maximum matching $M$ that contains the edge $xy$, for then $M \setminus \{xy\}$ will be a matching of $G'$, and consequently, $\nu(G') \geq |M \setminus \{xy\}| = |M|-1 = \nu(G)-1$, which is what we need. 

Now, fix any maximum matching $M$ of $G$. We may assume that $xy \notin M$, for otherwise we are done. It is clear that at least one of $x$, $y$ is $M$-saturated, for otherwise, $M \cup \{xy\}$ would be a matching of $G$ of size $|M|+1$, contrary to the maximality of $M$. Suppose first that exactly one of $x$ and $y$ is $M$-saturated; by symmetry, we may assume that $x$ is $M$-saturated, while $y$ is not. Then there exists some $x' \in V(G) \setminus \{x,y\}$ such that $xx' \in M$. But now $M' := (M \setminus \{xx'\}) \cup \{xy\}$ is a matching of size $|M| = \nu(G)$, and by construction, $xy \in M'$. We may now assume that both $x$ and $y$ are $M$-saturated. Since $xy \notin M$, we see that there exist distinct $x',y' \in V(G) \setminus \{x,y\}$ such that $xx',yy' \in M$. Since $N_G[x] = N_G[y] = C$, we see that $x',y' \in C$, and since $C$ is a clique, we see that $x'y' \in E(G)$. But now $(M \setminus \{xx',yy'\}) \cup \{xy,x'y'\}$ is a matching of $G$ of size $|M| = \nu(G)$, and by construction, $xy$ belongs to this new matching. This completes the proof of~(2).~$\blacklozenge$ 

\begin{adjustwidth}{1cm}{1cm} 
\begin{claim} \label{prop-delete-edge-claim-x-y-not-iso} Neither $x$ nor $y$ is an internal vertex of any induced path in $G$. Consequently, neither $x$ nor $y$ is an internal vertex of any isometric path in $G$.
\end{claim} 
\end{adjustwidth} 
\noindent 
{\em Proof of Claim~\ref{prop-delete-edge-claim-x-y-not-iso}.} Since all isometric paths are induced, it suffices to prove the first statement. Fix an induced path $P$ in $G$. Then any internal vertex of $P$ has two nonadjacent neighbors in $P$, and therefore in $G$ as well (because the path $P$ is induced). Since $x$ and $y$ are simplicial, it follows that neither $x$ nor $y$ is an internal vertex of $P$.~$\blacklozenge$

\begin{adjustwidth}{1cm}{1cm} 
\begin{claim} \label{prop-delete-edge-claim-isometric-in-G'} Any induced (resp.\ isometric) path of $G'$ is also an induced (resp.\ isometric) path in $G$. 
\end{claim} 
\end{adjustwidth} 
\noindent 
{\em Proof of Claim~\ref{prop-delete-edge-claim-isometric-in-G'}.} Since $G'$ is an induced subgraph of $G$, it is clear that all induced paths in $G'$ are also induced paths in $G$. Now, fix an isometric path $P$ of $G'$, and let $p$ and $q$ be its endpoints (possibly $p = q$). We must show that the path $P$ is isometric in $G$. Suppose otherwise, and fix an isometric path $Q$ between $p$ and $q$ in $G$ such that $\text{length}(Q) < \text{length}(P)$. Then at least one of $x$ and $y$ must be an internal vertex of $Q$, for otherwise, $Q$ would be a path in $G'$, contrary to the fact that the path $P$ is isometric in $G'$. But this contradicts Claim~\ref{prop-delete-edge-claim-x-y-not-iso}.~$\blacklozenge$

\begin{adjustwidth}{1cm}{1cm} 
\begin{claim} \label{prop-delete-edge-claim-xy-in-some-IPP-of-G} Some minimum IPP of $G$ contains the path $xy$. Similarly, some minimum IndPP of $G$ contains the path $xy$. 
\end{claim} 
\end{adjustwidth} 
\noindent 
{\em Proof of Claim~\ref{prop-delete-edge-claim-xy-in-some-IPP-of-G}.} We prove the claim for IPP; the proof for IndPP is completely analogous.\footnote{Indeed, as the reader can check, our argument remains valid if we simply replace all instances of ``IPP''  by ``IndPP,'' all instances of ``$\text{ipp}(G)$'' by ''$\text{indpp}(G)$,'' and all instances of ``isometric'' by ``induced.''} 

First of all, it is clear that $xy$ is, in fact, an isometric path in $G$. Now, fix a minimum IPP $\mathcal{I}$ of $G$, and fix $P_x,P_y \in \mathcal{I}$ such that $x \in V(P_x)$ and $y \in V(P_y)$. Since the paths $P_x$ and $P_y$ are isometric paths in~$G$, they are induced paths of $G$. Moreover, by Claim~\ref{prop-delete-edge-claim-x-y-not-iso}, $x$ is an endpoint of $P_x$, and $y$ is an endpoint of~$P_y$. 

Suppose first that $P_x = P_y$. Then $x$ and $y$ are the endpoints of $P_x = P_y$; since $x$ and $y$ are adjacent, and the path $P_x = P_y$ is isometric in $G$, we deduce that $P_x = P_y = xy$. Thus, $xy \in \mathcal{I}$, and we are done. 

From now on, we assume that $P_x \neq P_y$. Suppose that at least one of $P_x$ and $P_y$ has length zero. By symmetry, we may assume that $\text{length}(P_x) = 0$. Now, if $\text{length}(P_y) = 0$, then $\big(\mathcal{I} \setminus \{P_x,P_y\}\big) \cup \{xy\}$ is an IPP of $G$ of size $|\mathcal{I}|-1$, contrary to the minimality of $\mathcal{I}$. So, $\text{length}(P_y) \geq 1$. Set $P_y = yp_1\dots p_t$ ($t \geq 1$). But now $\big(\mathcal{I} \setminus \{P_x,P_y\}\big) \cup \{xy,p_1\dots p_t\}$ is an IPP of $G$ of size $|\mathcal{I}| = \text{ipp}(G)$, and by construction, this IPP contains~$xy$. 

It remains to consider the case when both $P_x$ and $P_y$ are of length at least one. Set $P_x = xx_1 \dots x_s$ and $P_y = yy_1\dots y_t$ ($s,t \geq 1$). Obviously, at most one of $P_x$ and $P_y$ contains the cut-vertex $v$; by symmetry, we may assume that $v \notin V(P_x)$. So, $P_x$ is in fact an induced path of the complete graph $G[C \setminus \{v\}]$, and we deduce that $\text{length}(P_x) = 1$, i.e.\ $P_x = xx_1$, with $x_1 \in C \setminus \{v\}$. But now $\big(\mathcal{I} \setminus \{P_x,P_y\}\big) \cup \{xy,x_1y_1\dots y_t\}$ is an IPP of $G$ of size $|\mathcal{I}| = \text{ipp}(G)$, and by construction, this IPP contains $xy$.~$\blacklozenge$ 

\begin{adjustwidth}{1cm}{1cm} 
\begin{claim} \label{prop-delete-edge-claim-ipp-G'-G-1} $\text{ipp}(G') = \text{ipp}(G)-1$ and $\text{indpp}(G') = \text{indpp}(G)-1$. 
\end{claim} 
\end{adjustwidth} 
\noindent 
{\em Proof of Claim~\ref{prop-delete-edge-claim-ipp-G'-G-1}.} We prove the first statement; the proof of the second statement is completely analogous.

We first show that $\text{ipp}(G') \geq \text{ipp}(G) - 1$. Fix a minimum IPP $\mathcal{I'}$ of $G'$. By Claim~\ref{prop-delete-edge-claim-isometric-in-G'}, every path in $\mathcal{I'}$ is isometric in $G$. But now $\mathcal{I'} \cup \{xy\}$ is an IPP of $G$. We deduce that $\text{ipp}(G) \leq |\mathcal{I'} \cup \{xy\}| = |\mathcal{I'}|+1 = \text{ipp}(G')+1$, and consequently, $\text{ipp}(G') \geq \text{ipp}(G)-1$. 

It remains to show that $\text{ipp}(G') \leq \text{ipp}(G)-1$. Using Claim~\ref{prop-delete-edge-claim-xy-in-some-IPP-of-G}, we fix a minimum IPP $\mathcal{I}$ of $G$ such that $xy \in \mathcal{I}$. But then $\mathcal{I} \setminus \{xy\}$ is an IPP of $G'$, and we deduce that $\text{ipp}(G') \leq |\mathcal{I} \setminus \{xy\}| = |\mathcal{I}|-1 = \text{ipp}(G)-1$.~$\blacklozenge$ 

\medskip 

By Claim~\ref{prop-delete-edge-claim-delete-edge}, we have that $|V(G')|-\nu(G') = |V(G)|-\nu(G)-1$. So, by Claim~\ref{prop-delete-edge-claim-ipp-G'-G-1}, we have the following: 
\begin{itemize} 
\item $\text{ipp}(G') = |V(G')|-\nu(G')$ if and only if $\text{ipp}(G) = |V(G)|-\nu(G)$;  
\item $\text{indpp}(G') = |V(G')|-\nu(G')$ if and only if $\text{indpp}(G) = |V(G)|-\nu(G)$. 
\end{itemize} 
Thus, both~(a) and~(b) hold, and we are done.  
\end{proof}

\begin{proposition} \label{prop-delete-odd-leaf-clique} Let $G$ be a graph that contains an odd leaf-clique $C$, and assume that $v$ is the unique cut-vertex of $G$ that belongs to $C$. Set $G' := G \setminus \big(C \setminus \{v\}\big)$. Then all the following hold: 
\begin{enumerate}[(a)] 
\item the blocks of $G$ are precisely $G[C]$ and the blocks of $G'$; 
\item $G$ is IPP-extremal if and only if $G'$ is IPP-extremal; 
\item $G$ is IndPP-extremal if and only if $G'$ is IndPP-extremal. 
\end{enumerate} 
\end{proposition} 
\begin{proof} 
Part~(a) is obvious, and we just need to prove~(b) and~(c). Since $C$ is an odd leaf-clique of $G$, we see that $|C| \geq 3$, and that $C \setminus \{v\}$ is a nonempty even clique. Set $k := \frac{|C \setminus \{v\}|}{2}$ and $C \setminus \{v\} = \{x_1,y_1,\dots,x_k,y_k\}$. So, we can obtain $G'$ from $G$ by iteratively deleting the sets $\{x_1,y_1\},\dots,\{x_k,y_k\}$ from $G$. Now~(b) and~(c) follow simply by applying Proposition~\ref{prop-delete-edge} $k$ times. 
\end{proof}

A graph is {\em chordal} if all induced cycles in it are triangles. The {\em diamond} is the graph on four vertices obtained from the complete graph $K_4$ by deleting one edge. A graph is {\em diamond-free} if none of its induced subgraphs is (isomorphic to) the diamond. 

\begin{theorem}\cite{bandelt1986distance} \label{thm-block-graph-chordal-diamond-free} Block graphs are precisely the diamond-free chordal graphs. 
\end{theorem}

\begin{proposition} \label{prop-odd-is-block} 
Let $G$ be a connected odd IPP-extremal graph. Then $G$ is a block graph. Moreover, if $G$ is biconnected, then $G$ is an odd complete graph. 
\end{proposition} 
\begin{proof} 
Clearly, biconnected block graphs are precisely the complete graphs. So, the second statement follows from the first (and from the assumption that $G$ is odd). Thus, it remains to show that $G$ is a block graph. In view of Theorem~\ref{thm-block-graph-chordal-diamond-free}, it suffices to show that $G$ is a diamond-free chordal graph. By Proposition~\ref{prop-basic-odd}, we have that $\text{ipp}(G) = \nu(G)+1$.

\begin{adjustwidth}{1cm}{1cm} 
\begin{claim} \label{prop-odd-is-block-claim-diamond-free} 
$G$ is diamond-free. 
\end{claim} 
\end{adjustwidth} 
\noindent 
{\em Proof of Claim~\ref{prop-odd-is-block-claim-diamond-free}.} Suppose otherwise, and fix an induced diamond $D$ in $G$, where $V(D) = \{a,b,c,d\}$ and $E(D) = \{ab,bc,cd,da,bd\}$. Using Proposition~\ref{prop-basic-odd}, we fix a perfect matching $M$ of $G \setminus d$. Then $M$ is a maximum matching in $G$ such that $d$ is the unique $M$-unsaturated vertex of $G$, and in particular, $a$, $b$, and $c$ are all $M$-saturated. Clearly, $M \cup \{d\}$ is an IPP of $G$. Since $d$ is $M$-unsaturated, we see that $E(D) \cap M \subseteq \{ab,bc\}$. But since edges $ab$ and $bc$ share an endpoint, the matching $M$ contains at most one of them. So, by symmetry, we may assume that either $E(D) \cap M = \{ab\}$ or $E(D) \cap M = \emptyset$. In each case, we will derive a contradiction by exhibiting an IPP of $G$ of size $\nu(G)$, contrary to the fact that $\text{ipp}(G) = \nu(G)+1$. 

Suppose first that $E(D) \cap M = \{ab\}$. Then, since $c$ is $M$-saturated, there exists some $c' \in V(G) \setminus V(D)$ such that $cc' \in M$. Since $d$ is $M$-unsaturated and adjacent to $c$, Proposition~\ref{prop-no-unsaturated-mixed-on-matching-edge} guarantees that $dc' \in E(G)$. But then $\big(M \setminus \{ab,cc'\}\big) \cup \{abc,dc'\}$ is an IPP of $G$ of size $|M| = \nu(G)$, a contradiction. 

Suppose now that $E(D) \cap M = \emptyset$. Then, since $a$, $b$, and $c$ are $M$-saturated, there exist pairwise distinct vertices $a',b',c' \in V(G) \setminus V(D)$ such that $aa', bb', cc' \in M$. Since $d$ is $M$-unsaturated and adjacent to $a$, $b$, and $c$, Proposition~\ref{prop-no-unsaturated-mixed-on-matching-edge} guarantees that $da', db', dc' \in E(G)$. If $a'b' \notin E(G)$, then $\big(M \setminus \{aa',bb'\}\big) \cup \{a'db', ab\}$ is an IPP of $G$ of size $|M| = \nu(G)$, a contradiction. So, $a'b' \in E(G)$. But now $\big(M \setminus \{aa',bb',cc'\}\big) \cup \{a'b', abc, dc'\}$ is an IPP of $G$ of size $|M| = \nu(G)$, again a contradiction.~$\blacklozenge$ 

\medskip 

It now remains to show that $G$ is chordal. Suppose otherwise, and fix an induced cycle $C = c_0c_1\dots c_{k-1}c_0$ (with $k \geq 4$ and with indices taken modulo $k$) of $G$.

\begin{adjustwidth}{1cm}{1cm} 
\begin{claim} \label{prop-odd-is-block-claim-diamond-free-claim-M-C} For every maximum matching $M$ of $G$, either all vertices of $C$ are $M$-saturated or $E(C) \cap M = \emptyset$. 
\end{claim} 
\end{adjustwidth}
\noindent 
{\em Proof of Claim~\ref{prop-odd-is-block-claim-diamond-free-claim-M-C}.} Suppose otherwise, and fix a maximum matching $M$ of $G$ such that for some $i,j \in \{0,1,\dots,k-1\}$, $c_i$ is $M$-unsaturated and $c_jc_{j+1} \in M$. Clearly, $i \notin \{j,j+1\}$. We may assume that $M$, $c_i$, and $c_jc_{j+1}$ were chosen so that the path $P := c_ic_{i+1}\dots c_j$ is as short as possible. By symmetry, we may assume that $i = 0$ (and consequently, $j \in \{1,\dots,k-2\}$). Proposition~\ref{prop-no-unsaturated-mixed-on-matching-edge} now guarantees that $c_0$ is the unique $M$-unsaturated vertex of $G$, and in particular, $c_1$ is $M$-saturated. So, there exists some $c_1' \in V(G) \setminus \{c_1\}$ such that $c_1c_1' \in M$. Since $c_0$ is $M$-unsaturated and adjacent to $c_1$, Proposition~\ref{prop-no-unsaturated-mixed-on-matching-edge} guarantees that $c_0c_1' \in E(G)$. Now $\{c_0,c_1,c_1'\}$ induces a triangle in $G$; since $C = c_0c_1\dots c_{k-1}c_0$ is an induced cycle of $G$ of length $k \geq 4$, we see that $c_1' \notin V(C)$, and in particular, $j \neq 1$. But now $M_1 := \big(M \setminus \{c_1c_1'\}\big) \cup \{c_0c_1'\}$ is a maximum matching of $G$, and $c_1$ is $M_1$-unsaturated. Meanwhile, we have that $c_jc_{j+1} \in M_1$, and the path $P' := c_1c_2\dots c_j$ is shorter than the path $P = c_0c_1\dots c_j$, contrary to the minimality of $P$.~$\blacklozenge$

\medskip 

Using Proposition~\ref{prop-basic-odd}, we fix a perfect matching $M$ of $G \setminus c_0$. By Claim~\ref{prop-odd-is-block-claim-diamond-free-claim-M-C}, we have that $E(C) \cap M = \emptyset$. We now deduce that there exist pairwise distinct vertices $c_1',\dots,c_{k-1}' \in V(G) \setminus V(C)$ such that $c_1c_1',\dots,c_{k-1}c_{k-1}' \in M$. 

\begin{adjustwidth}{1cm}{1cm} 
\begin{claim} \label{prop-odd-is-block-claim-diamond-free-claim-triangles} $c_0c_1',c_0c_{k-1}' \in E(G)$ and $c_{k-2}c_{k-1}' \notin E(G)$. 
\end{claim} 
\end{adjustwidth} 
\noindent
{\em Proof of Claim~\ref{prop-odd-is-block-claim-diamond-free-claim-triangles}.} Since the $M$-unsaturated vertex $c_0$ is adjacent to the endpoint $c_1$ of the matching edge $c_1c_1' \in M$, Proposition~\ref{prop-no-unsaturated-mixed-on-matching-edge} guarantees that $c_0c_1' \in E(G)$. Analogously, $c_0c_{k-1}' \in E(G)$. Finally, we have that $c_{k-2}c_{k-1}' \notin E(G)$, for otherwise, vertices $c_0,c_{k-2},c_{k-1},c_{k-1}'$ would induce a diamond in $G$, contrary to the fact that, by Claim~\ref{prop-odd-is-block-claim-diamond-free}, $G$ is diamond-free.~$\blacklozenge$ 

\medskip 

By Claim~\ref{prop-odd-is-block-claim-diamond-free-claim-triangles}, we have that $c_{k-2}c_{k-1}' \notin E(G)$. Now, fix the smallest index $i \in \{1,\dots,k-2\}$ such that $c_ic_{i+1}' \notin E(G)$. By the minimality of $i$, and by Claim~\ref{prop-odd-is-block-claim-diamond-free-claim-triangles}, we have that $c_jc_{j+1}' \in E(G)$ for all $j \in \{0,\dots,i-1\}$. But now 
\begin{displaymath} 
\begin{array}{rcl} 
\mathcal{I} & := & \{c_jc_{j+1}' \mid 0 \leq j \leq i-1\} \cup \big(M \setminus \{c_jc_j' \mid 1 \leq j \leq i+1\}\big) \cup \{c_ic_{i+1}c_{i+1}'\} 
\end{array} 
\end{displaymath} 
is an IPP of $G$ of size $\nu(G)$, contrary to the fact that $\text{ipp}(G) = \nu(G)+1$.
\end{proof}  

\begin{proposition} \label{prop-one-even-block-implies-extremal} Let $G$ be a connected graph that has at most one block that is not an odd complete graph. Assume furthermore that if $B$ is a block of $G$ that is not an odd complete graph, then $B$ is even and IPP-extremal. Then $G$ is IPP-extremal. 
\end{proposition} 
\begin{proof} 
We may assume inductively that for any connected graph $G'$ such that $|V(G')| < |V(G)|$, if the following hold: 
\begin{itemize} 
\item $G'$ has at most one block that is not an odd complete graph, and 
\item if $B'$ is a block of $G'$ that is not an odd complete graph, then $B'$ is even and IPP-extremal, 
\end{itemize} 
then $G'$ is IPP-extremal. 

Suppose first that $G$ is biconnected, so that $G$ has exactly one block, namely $G$ itself. If $G$ is even, then the result is immediate. On the other hand, if $G$ is odd, then by assumption, $G$ is an odd complete graph, and we are done (because clearly, all complete graphs are IPP-extremal).  %Proposition~\ref{prop-odd-is-block}  guarantees that $G$ is an odd complete graph, and we are done. 

Suppose now that $G$ is not biconnected. Then $G$ has at least two leaf-blocks, and by hypothesis, at least one of them is an odd complete graph. Let $C$ be the vertex set of this block, and let $v$ be the unique cut-vertex of $G$ that belongs to $C$. Set $G' := G \setminus \big(C \setminus \{v\}\big)$. Clearly, $G'$ is connected, and by Proposition~\ref{prop-delete-odd-leaf-clique}(a), the blocks of $G$ are precisely $G[C]$ and the blocks of $G'$. So, by the induction hypothesis, $G'$ is IPP-extremal. But now Proposition~\ref{prop-delete-odd-leaf-clique}(b) guarantees that $G$ is IPP-extremal. 
\end{proof} 

\begin{proposition} \label{prop-one-even-block-implies-extremal-IndPP} Let $G$ be a connected graph that has at most one block that is not an odd complete graph. Assume furthermore that if $B$ is a block of $G$ that is not an odd complete graph, then $B$ is even and IndPP-extremal. Then $G$ is IndPP-extremal. 
\end{proposition} 
\begin{proof} 
The proof is completely analogous to that of Proposition~\ref{prop-one-even-block-implies-extremal}: we simply replace all instances of ``IPP'' by ``IndPP,'' and we use Proposition~\ref{prop-delete-odd-leaf-clique}(c) instead of Proposition~\ref{prop-delete-odd-leaf-clique}(b).%\footnote{Note that Proposition~\ref{prop-odd-is-block} (used in the proof of Proposition~\ref{prop-one-even-block-implies-extremal}) does in fact apply to IndPP-extremal graphs, since by Proposition~\ref{prop-indpp-implies-ipp}, IndPP-extremal graphs are IPP-extremal.}
\end{proof} 

Proposition~\ref{prop-one-even-block-implies-extremal} yields the backward direction of Theorem~\ref{thm-main}, and Proposition~\ref{prop-one-even-block-implies-extremal-IndPP} yields the backward direction of Theorem~\ref{thm-main-indpp}. It remains to prove the forward direction of the two theorems. 

\medskip 

%A {\em short independent path partition} (or {\em SIndPP}) of a graph $G$ is an IndPP $\mathcal{I}$ of $G$ such that all paths in $\mathcal{I}$ are of length at most two (i.e.\ have at most three vertices). Since all induced paths of length at least two are isometric, it clear that any SIndPP is an IPP. 

Given a graph $G$ and a vertex $v \in V(G)$, we say that an IPP $\mathcal{I}$ of $G$ is {\em $v$-extendable} if $|\mathcal{I}| \leq \nu(G)$ and $v$ is an endpoint of the unique path in $\mathcal{I}$ that contains $v$ (possibly, this path is a one-vertex path whose only vertex is $v$). 

\begin{proposition} \label{prop-v-extend} Let $G$ be a biconnected graph, and let $v \in V(G)$. Assume that $G \setminus v$ has a perfect matching, and that $\text{ipp}(G) \leq \nu(G)$. Then $G$ admits a $v$-extendable IPP. 
\end{proposition} 
\begin{proof} 
We assume toward a contradiction that $G$ does not admit a $v$-extendable IPP. 

\begin{adjustwidth}{1cm}{1cm} 
\begin{claim} \label{prop-v-extend-claim-v-not-mixed-on-matching-edge} The vertex $v$ is not mixed on any edge of any perfect matching of $G \setminus v$. 
\end{claim} 
\noindent 
\end{adjustwidth} 
{\em Proof of Claim~\ref{prop-v-extend-claim-v-not-mixed-on-matching-edge}.} Fix a perfect matching $M$ of $G \setminus v$, and assume toward a contradiction that $v$ is mixed on some edge $xy \in M$. By symmetry, we may assume that $v$ is adjacent to $x$ and nonadjacent to~$y$. But then $\mathcal{I} := \big(M \setminus \{xy\}\big) \cup \{vxy\}$ is a $v$-extendable IPP of $G$, a contradiction.~$\blacklozenge$

\begin{adjustwidth}{1cm}{1cm} 
\begin{claim} \label{prop-v-extend-claim-NGv-not-mixed-on-M} 
For every edge $xx'$ that belongs to some perfect matching of $G \setminus v$, no vertex in $N_G(v) \setminus \{x,x'\}$ is mixed on $xx'$. 
\end{claim} 
\end{adjustwidth} 
{\em Proof of Claim~\ref{prop-v-extend-claim-NGv-not-mixed-on-M}.} Fix a perfect matching $M$ of $G \setminus v$, and fix an edge $xx' \in M$. Suppose toward a contradiction that some vertex $u \in N_G(v) \setminus \{x,x'\}$ is mixed on $xx'$. By symmetry, we may assume that $u$ is adjacent to $x$ and nonadjacent to $x'$. Since $M$ is a perfect matching of $G \setminus v$, we know that there exists some vertex $u' \in V(G) \setminus \{v,u,x,x'\}$ such that $uu' \in M$. Since $u \in N_G(v)$, Claim~\ref{prop-v-extend-claim-v-not-mixed-on-matching-edge} guarantees that $u' \in N_G(v)$. But now $\mathcal{I} := \big(M \setminus \{xx',uu'\}\big) \cup \{uxx',vu'\}$ is a $v$-extendable IPP of $G$, a contradiction.~$\blacklozenge$ 

\begin{adjustwidth}{1cm}{1cm} 
\begin{claim} \label{prop-v-extend-claim-M-in-NGv-not-mixed} For all perfect matchings $M$ of $G \setminus v$, and all distinct edges $xx',yy' \in M$ such that $x,x',y,y' \in N_G(v)$, we have that $xx'$ and $yy'$ are either complete or anticomplete to each other. 
\end{claim} 
\end{adjustwidth} 
{\em Proof of Claim~\ref{prop-v-extend-claim-M-in-NGv-not-mixed}.} Fix a perfect matching $M$ of $G \setminus v$, and fix distinct edges $xx',yy' \in M$ such that $x,x',y,y' \in N_G(v)$. We may assume that $xx'$ is not anticomplete to $yy'$, for otherwise we are done. By symmetry, we may assume that $x$ is adjacent to $y$. Since $x \in N_G(v)$ and $yy' \in M$, Claim~\ref{prop-v-extend-claim-NGv-not-mixed-on-M} guarantees that $x$ is complete to $yy'$. In particular, $y' \in N_G(v)$ is adjacent to $x$; so, since $xx' \in M$, Claim~\ref{prop-v-extend-claim-NGv-not-mixed-on-M} guarantees that $y'$ is complete to $xx'$. But now $x' \in N_G(v)$ is adjacent to $y'$, and so since $yy' \in M$, Claim~\ref{prop-v-extend-claim-NGv-not-mixed-on-M} guarantees that $x'$ is complete to $yy'$. We have now shown that $xx'$ is complete to $yy'$, and we are done.~$\blacklozenge$

\begin{adjustwidth}{1cm}{1cm} 
\begin{claim} \label{prop-v-extend-claim-NGv-cliques} 
$N_G(v)$ can be partitioned into two or more nonempty cliques, pairwise anticomplete to each other. 
\end{claim} 
\end{adjustwidth} 
{\em Proof of Claim~\ref{prop-v-extend-claim-NGv-cliques}.} It suffices to show that $v$ is not simplicial in $G$, and that $G[N_G(v)]$ contains no induced three-vertex path.\footnote{Let us explain why this is enough. Suppose that we have proven that $v$ is not simplicial in $G$ (i.e.\ $N_G(v)$ is not a clique, and in particular, $N_G(v) \neq \emptyset$), and that $G[N_G(v)]$ contains no induced three-vertex path. Let $C$ be a connected component of $G[N_G(v)]$, and suppose that $C$ is not a complete graph. Then $C$ contains two nonadjacent vertices, call them $c_1$ and $c_2$. Fix an induced path $P$ between $c_1$ and $c_2$ in $C$. This path contains at least three vertices, and so some induced subpath $P'$ of $P$ is a three-vertex path. Clearly, the three-vertex $P'$ is an induced path in $G[N_G(v)]$, a contradiction. Thus, all connected components of $G[N_G(v)]$ are complete, and since $v$ is not simplicial in $G$, it follows that $G[N_G(v)]$ has at least two connected components. In other words, $N_G(v)$ can be partitioned into two or more nonempty cliques, pairwise anticomplete to each other.} 

First, suppose toward a contradiction that $v$ is simplicial in $G$. Fix any minimum IPP $\mathcal{I}$ of $G$; by hypothesis, we have that $|\mathcal{I}| = \text{ipp}(G) \leq \nu(G)$. Fix the (unique) path $P \in \mathcal{I}$ such that $v \in V(P)$. Clearly, every internal vertex of $P$ has two nonadjacent neighbors in $P$; since the path $P$ is isometric and therefore induced in $G$, it follows that no internal vertex of $P$ is simplicial in $G$. So, $v$ is not an internal vertex of $P$, and it must therefore be an endpoint of $P$. But now $\mathcal{I}$ is a $v$-extendable IPP of $G$, a contradiction. This proves that $v$ is not simplicial in $G$. 

Now, suppose toward a contradiction that $G[N_G(v)]$ contains an induced three-vertex path $xyz$. Fix a perfect matching $M$ of $G \setminus v$ (such an $M$ exists by hypothesis). Since $z \in N_G(v)$ is mixed on $xy$, Claim~\ref{prop-v-extend-claim-NGv-not-mixed-on-M} guarantees that $xy \notin M$. Similarly, since $x \in N_G(v)$ is mixed on $yz$, Claim~\ref{prop-v-extend-claim-NGv-not-mixed-on-M} guarantees that $yz \notin M$. Since $M$ is a perfect matching of $G \setminus v$, it follows that there exist pairwise distinct vertices $x',y',z' \in V(G) \setminus \{v,x,y,z\}$ such that $xx',yy',zz' \in M$. Since $x,y,z \in N_G(v)$, Claim~\ref{prop-v-extend-claim-v-not-mixed-on-matching-edge} guarantees that $x',y',z' \in N_G(v)$. 

We now have that $x,x',y,y',z,z' \in N_G(v)$, that $xx',yy',zz' \in M$, and that $xyz$ is an induced path in $G$. So, by Claim~\ref{prop-v-extend-claim-M-in-NGv-not-mixed}, we have that $\{y,y'\}$ is complete to $\{x,x',z,z'\}$, and that $\{x,x'\}$ and $\{z,z'\}$ are anticomplete to each other. In particular, $x'y'z'$ is an induced path in $G$. But we now see that $\mathcal{I} := \big(M \setminus \{xx',yy',zz'\}\big) \cup \{xyz,x'y'z',v\}$ is a $v$-extendable IPP of $G$, a contradiction. This proves that $G[N_G(v)]$ contains no induced three-vertex path, and we are done.~$\blacklozenge$ 

\medskip 

Using Claim~\ref{prop-v-extend-claim-NGv-cliques}, we fix a partition $(C_1,\dots,C_k)$, with $k \geq 2$, of $N_G(v)$ into nonempty cliques, pairwise anticomplete to each other. 

\begin{adjustwidth}{1cm}{1cm} 
\begin{claim} \label{prop-v-extend-claim-matching-edge-in-Ci} For every perfect matching $M$ of $G \setminus v$, and every edge $xx' \in M$, either $x,x' \in V(G) \setminus N_G[v]$, or there exists some $i \in \{1,\dots,k\}$ such that $x,x' \in C_i$. 
\end{claim}  
\end{adjustwidth} 
{\em Proof of Claim~\ref{prop-v-extend-claim-matching-edge-in-Ci}.} Fix a perfect matching $M$ of $G \setminus v$, and fix an edge $xx' \in M$. By Claim~\ref{prop-v-extend-claim-v-not-mixed-on-matching-edge}, either $x,x' \in N_G(v)$ or $x,x' \in V(G) \setminus N_G[v]$. We may assume that $x,x' \in N_G(v)$, for otherwise we are done. The result now follows from the fact that $C_1,\dots,C_k$ form a partition of $N_G(v)$ into cliques, pairwise anticomplete to each other.~$\blacklozenge$

\begin{adjustwidth}{1cm}{1cm} 
\begin{claim} \label{prop-v-extend-claim-max-one-neighbor-in-C_i} No vertex in $V(G) \setminus N_G[v]$ has more than one neighbor in any one of the cliques $C_1,\dots,C_k$. 
\end{claim} 
\end{adjustwidth} 
\noindent 
{\em Proof of Claim~\ref{prop-v-extend-claim-max-one-neighbor-in-C_i}.} Suppose otherwise. By symmetry, we may assume that some vertex $u \in V(G) \setminus N_G[v]$ has at least two neighbors, call them $x$ and $y$, in $C_1$. 

We first construct a perfect matching $M$ of $G \setminus v$ such that $xy \in M$. Fix any perfect matching $M_0$ of $G \setminus v$ (such an $M_0$ exists by hypothesis). If $xy \in M_0$, then we simply set $M := M_0$. So, suppose that $xy \notin M_0$. Then there exist distinct vertices $x',y' \in V(G) \setminus \{v,x,y\}$ such that $xx',yy' \in M_0$. Since $x,y \in C_1$, Claim~\ref{prop-v-extend-claim-matching-edge-in-Ci} guarantees that $x',y' \in C_1$. But $C_1$ is a clique, and so $x'y' \in E(G)$. We now set $M := \big(M_0 \setminus \{xx',yy'\}\big) \cup \{xy,x'y'\}$, and we observe that $M$ is a perfect matching of $G \setminus v$ such that $xy \in M$. 

Now, since $M$ is a perfect matching of $G \setminus v$, there exists some $u' \in V(G) \setminus \{v,u,x,y\}$ such that $uu' \in M$. Since $u \in V(G) \setminus N_G[v]$, Claim~\ref{prop-v-extend-claim-v-not-mixed-on-matching-edge} guarantees that $u' \in V(G) \setminus N_G[v]$. Since $x,y \in N_G(v)$ are both adjacent to $u$, and since $uu' \in M$, Claim~\ref{prop-v-extend-claim-NGv-not-mixed-on-M} guarantees that $\{x,y\}$ is complete to $\{u,u'\}$. But now $\mathcal{I} := \big(M \setminus \{xy,uu'\}\big) \cup \{vxu,yu'\}$ is a $v$-extendable IPP of $G$, a contradiction.~$\blacklozenge$ 

\medskip 

Now, since $G$ is biconnected, we see that $G \setminus v$ is connected. So, there exists an induced path $Q := q_0q_1\dots q_tq_{t+1}$ ($t \geq 1$) in $G \setminus v$ such that $q_0 \in C_1$ and $q_{t+1} \in C_2 \cup \dots \cup C_k$, whereas $q_1,\dots,q_t$ (the internal vertices of $Q$) belong to $V(G) \setminus N_G[v]$. By symmetry, we may assume that $q_{t+1} \in C_2$. Further, fix a perfect matching $M$ of $G \setminus v$ (such an $M$ exists by hypothesis), and we fix vertices $q_0',q_1',\dots,q_t',q_{t+1}' \in V(G) \setminus \{v\}$ such that $q_0q_0',q_1q_1',\dots,q_tq_t',q_{t+1}q_{t+1}' \in M$.

\begin{adjustwidth}{1cm}{1cm} 
\begin{claim} \label{prop-v-extend-claim-Q-spikes} All the following hold: 
\begin{enumerate}[(1)] 
\item vertices $q_0,q_1,\dots,q_t,q_{t+1},q_0',q_1',\dots,q_t',q_{t+1}'$ are pairwise distinct (and in particular, no edge of $Q$ belongs to $M$); 
\item for all $i \in \{0,\dots,t-1\}$, we have that $q_iq_{i+1}' \in E(G)$; 
\item $q_tq_{t+1}' \notin E(G)$. 
\end{enumerate} 
\end{claim} 
\end{adjustwidth}
\noindent 
{\em Proof of Claim~\ref{prop-v-extend-claim-Q-spikes}.} It is clear that $q_0,q_1,\dots,q_t,q_{t+1}$ are pairwise distinct, and it is also clear that $q_0',q_1',\dots,q_t',q_{t+1}'$ are pairwise distinct. Thus, to prove~(1), it is enough to show that $q_0',q_1',\dots,q_t',q_{t+1}' \notin V(Q)$. 

First, since $q_0 \in C_1$, $q_{t+1} \in C_2$, and $q_0q_0',q_{t+1}q_{t+1}' \in M$, Claim~\ref{prop-v-extend-claim-matching-edge-in-Ci} guarantees that $q_0' \in C_1$ and $q_{t+1}' \in C_2$. In particular, $q_0',q_{t+1}' \notin V(Q)$. Furthermore, since $q_1,\dots,q_t \in V(G) \setminus N_G[v]$, Claim~\ref{prop-v-extend-claim-v-not-mixed-on-matching-edge} guarantees that $q_1',\dots,q_t' \in V(G) \setminus N_G[v]$. Further, since $q_{t+1},q_{t+1}' \in C_2$, and since $q_t \in V(G) \setminus N_G[v]$ is adjacent to $q_{t+1}$, Claim~\ref{prop-v-extend-claim-max-one-neighbor-in-C_i} guarantees that $q_t$ is nonadjacent to $q_{t+1}'$; this proves~(3). It remains to prove~(1) and~(2). 

Since $q_0 \in N_G(v)$ is adjacent to $q_1$, and since $q_1q_1' \in M$, Claim~\ref{prop-v-extend-claim-NGv-not-mixed-on-M} guarantees that $q_0q_1' \in E(G)$. Consequently, $q_1' \notin V(Q)$ (because the path $Q$ is induced). We have now shown that $q_0',q_1' \notin V(Q)$ and $q_0q_1' \in E(G)$. Fix the largest index $s \in \{0,\dots,t-1\}$ such that $q_0',\dots,q_{s+1}' \notin V(Q)$, and such that for all indices $i \in \{0,\dots,s\}$, we have that $q_iq_{i+1}' \in E(G)$. We have already shown that $q_{t+1}' \notin V(Q)$, and so if $s = t-1$, then~(1) and~(2) hold, and we are done. We may therefore assume that $s \leq t-2$. 

Since $q_{s+1}' \notin V(Q)$, we know that $q_{s+2}' \neq q_{s+1}$.\footnote{Indeed, edges $q_{s+1}q_{s+1}'$ and $q_{s+2}q_{s+2}'$ both belong to the matching $M$, and so either these two edges are equal, or they have no endpoints in common. So, if we had $q_{s+2}' = q_{s+1}$, then it would follow that $q_{s+2} = q_{s+1}'$, contrary to the fact that $q_{s+2} \in V(Q)$ and $q_{s+1}' \notin V(Q)$.} But now $q_{s+1}q_{s+2}' \in E(G)$, for otherwise, 
\begin{displaymath} 
\begin{array}{rcl} 
\mathcal{I} & := & \big(M \setminus \{q_0q_0',q_1q_1',\dots,q_{s+2}q_{s+2}'\}\big) \cup \{vq_0',q_0q_1',\dots,q_sq_{s+1}'\} \cup \{q_{s+1}q_{s+2}q_{s+2}'\} 
\end{array} 
\end{displaymath} 
would be a $v$-extendable IPP of $G$, a contradiction. Since the path $Q$ is induced in $G$, this further implies that $q_{s+2}' \notin V(Q)$. We have now derived a contradiction to the maximality of the index $s$.~$\blacklozenge$

\medskip 

Using Claim~\ref{prop-v-extend-claim-Q-spikes}, we now see that 
\begin{displaymath} 
\begin{array}{rcl} 
\mathcal{I} & := & \big(M \setminus \{q_0q_0',q_1q_1',\dots,q_tq_t',q_{t+1}q_{t+1}'\}\big) \cup \{vq_0',q_0q_1',\dots,q_{t-1}q_t',q_tq_{t+1}q_{t+1}'\} 
\end{array} 
\end{displaymath} 
is a $v$-extendable IPP of $G$, a contradiction. This completes the argument. 
\end{proof}

\begin{proposition} \label{prop-extremal-implies-one-even-block} Let $G$ be a connected IPP-extremal graph. Then $G$ has at most one block that is not an odd complete graph. Moreover, if $G$ has exactly one block that is not an odd complete graph, then that block is even. %Moreover, if $G$ has a block $B$ that is not an odd complete graph, then the block $B$ is even and satisfies $\text{ipp}(B) = |V(B)|-\nu(B)$. 
\end{proposition} 
\begin{proof} 
We may assume inductively that the proposition is true for graphs on fewer than $|V(G)|$ many vertices. More precisely, we assume inductively that for all connected IPP-extremal graphs $G'$ such that $|V(G')| < |V(G)|$, the graph $G'$ has at most one block that is not an odd complete graph, and moreover, if $G'$ has exactly one block that is not an odd complete graph, then that block is even. 

%if $\text{ipp}(G') = |V(G')|-\nu(G')$, then the following hold: 
%\begin{itemize} 
%\item $G'$ has at most one block that is not an odd complete graph, and 
%\item if $G'$ has a block $B'$ that is not an odd complete graph, then the block $B'$ is even and satisfies $\text{ipp}(B') = |V(B')|-\nu(B')$. 
%\end{itemize} 

Suppose first that $G$ is biconnected, so that $G$ has exactly one block, namely $G$ itself. If $G$ is even, then the result is immediate, whereas if $G$ is odd, then Proposition~\ref{prop-odd-is-block} guarantees that $G$ is an odd complete graph, and we are done. 

From now on, we assume that $G$ is not biconnected. Therefore, $G$ has at least two leaf-blocks.

\begin{adjustwidth}{1cm}{1cm} 
\begin{claim} \label{prop-extremal-implies-one-even-block-claim-one-even-block-leaf}
$G$ has at most one even leaf-block. 
\end{claim} 
\end{adjustwidth} 
{\em Proof of Claim~\ref{prop-extremal-implies-one-even-block-claim-one-even-block-leaf}.} Suppose otherwise, and fix distinct even leaf-blocks $B_1$ and $B_2$ of $G$. Further, for each index $i \in \{1,2\}$, let $v_i$ be the (unique) cut-vertex of $G$ that belongs to $B_i$ (it is possible that $v_1 = v_2$). 

Suppose first that $G$ is odd. Then Proposition~\ref{prop-basic-odd} guarantees that $G \setminus v_1$ admits a perfect matching. But this is impossible because $G \setminus v_1$ has at least one odd connected component, namely $B_1 \setminus v_1$. 

Thus, $G$ is even. Then Proposition~\ref{prop-basic-even} guarantees that $G$ admits a perfect matching, call it $M$. For each $i \in \{1,2\}$, we fix $b_i \in V(G) \setminus \{v_i\}$ such that $b_iv_i \in M$, and we observe that $b_i \in V(B_i) \setminus \{v_i\}$, for otherwise, $B_i \setminus v_i$ would be an odd connected component of $G \setminus \{b_i,v_i\}$, contrary to the fact that $G \setminus \{b_i,v_i\}$ has a perfect matching (namely $M \setminus \{b_iv_i\}$). In particular, matching edges $b_1v_1$ and $b_2v_2$ are distinct (and consequently have no common endpoints). Now, let $P$ be a shortest path between $b_1$ and $b_2$ in $G$. Clearly, the path $P$ contains the matching edges $b_1v_1$ and $b_2v_2$. Set $M' := \{uv \in M \mid u,v \notin V(P)\}$, and note that $M' \subseteq M \setminus \{b_1v_1,b_2v_2\}$. But now 
\begin{displaymath} 
\begin{array}{rcl} 
\mathcal{I} & := & \{P\} \cup M' \cup \big\{v \in V(G) \setminus V(P) \mid \text{$\exists u \in V(P)$ s.t.\ $uv \in M$}\big\}
\end{array} 
\end{displaymath} 
is an IPP of $G$, and consequently, $\text{ipp}(G) \leq |\mathcal{I}| \leq |M|-1 = \nu(G)-1$, contrary to Proposition~\ref{prop-basic-even}.~$\blacklozenge$  

\medskip 

Since $G$ has at least two leaf-blocks, and since at most one of them is even (by Claim~\ref{prop-extremal-implies-one-even-block-claim-one-even-block-leaf}), we see that $G$ has an odd leaf-block, call it $B$. Let $v$ be the (unique) cut-vertex of $G$ that belongs to $B$. Now, if $G$ is even, then set $U := \emptyset$; and if $G$ is odd, then let $U$ be a singleton whose only member is some vertex of $V(G) \setminus V(B)$. Propositions~\ref{prop-basic-even} and~\ref{prop-basic-odd} now guarantee that $G_U := G \setminus U$ has a perfect matching, call it $M$, and that $M \cup U$ is a minimum IPP of $G$, i.e.\ $\text{ipp}(G) = |M|+|U|$. By construction, we have that $v \notin U$, and so $v$ is $M$-saturated; fix $v^* \in V(G) \setminus \{v\}$ such that $vv^* \in M$. Set $M_B := M \cap E(B)$ and $M_B' := M \setminus E(B)$. 

\begin{adjustwidth}{1cm}{1cm} 
\begin{claim} \label{prop-extremal-implies-one-even-block-claim-MB} We have that $v^* \notin V(B)$ and $vv^* \in M_B'$. Moreover, $M_B$ is a perfect matching of $B \setminus v$, and consequently, $\nu(B) = |M_B|$. 
\end{claim} 
\end{adjustwidth} 
{\em Proof of Claim~\ref{prop-extremal-implies-one-even-block-claim-MB}.} If we had that $v^* \in V(B)$, then some connected component $C$ of $B \setminus \{v,v^*\}$ would be odd, and clearly, $C$ would be an odd connected component of $G_U \setminus \{v,v^*\}$ as well, contrary to the fact that $G_U \setminus \{v,v^*\}$ has a perfect matching, namely $M \setminus \{vv^*\}$. Thus, $v^* \notin V(B)$. It follows that, $vv^* \notin M_B$, and consequently, $vv^* \in M_B'$. But now $B \setminus v$ is a connected component of $G_U \setminus \{v,v^*\}$, and $M \setminus \{vv^*\}$ is a perfect matching of $G_U \setminus \{v,v^*\}$. So, $\big(M \setminus \{vv^*\}\big) \cap E(B)$ is a perfect matching of $B \setminus v$, and clearly, $\big(M \setminus \{vv^*\}\big) \cap E(B) = \big(M \cap E(B)\big) \setminus \{vv^*\} = M_B \setminus \{vv^*\} = M_B$.~$\blacklozenge$ 

\begin{adjustwidth}{1cm}{1cm} 
\begin{claim} \label{prop-extremal-implies-one-even-block-claim-B-extremal} $B$ is IPP-extremal. 
\end{claim}  
\end{adjustwidth} 
{\em Proof of Claim~\ref{prop-extremal-implies-one-even-block-claim-B-extremal}.} 
Suppose otherwise, so that (by Proposition~\ref{prop-upper-bound-trivial}) we have that $\text{ipp}(B) \leq |V(B)|-\nu(B)-1$. Since $B \setminus v$ has a perfect matching (by Claim~\ref{prop-extremal-implies-one-even-block-claim-MB}), we see that $|V(B)| = 2\nu(B)+1$, and we deduce that $\text{ipp}(B) \leq \nu(B)$. 
% We now know that $B$ is biconnected (because it is a leaf-block of $G$), that $B \setminus v$ admits a perfect matching (by Claim~\ref{prop-extremal-implies-one-even-block-claim-MB}), and that $\text{ipp}(B) \leq \nu(B)$. So, 
Hence, by Proposition~\ref{prop-v-extend}, $B$ admits a $v$-extendable IPP $\mathcal{I}_B$. By the definition of a $v$-extendable IPP, and by Claim~\ref{prop-extremal-implies-one-even-block-claim-MB}, we have that $|\mathcal{I}_B| \leq \nu(B) = |M_B|$. Further, fix $P \in \mathcal{I}_B$ such that $v \in V(P)$. Since $\mathcal{I}_B$ is $v$-extendable, we know that $v$ is an endpoint of $P$; set $P = vp_1\dots p_t$ ($t \geq 0$). By Claim~\ref{prop-extremal-implies-one-even-block-claim-MB}, we have that $v^* \notin V(B)$ and $vv^* \in M_B'$, and we deduce that 
\begin{displaymath} 
\begin{array}{rcl} 
\mathcal{I} & := & \{v^*vp_1\dots p_t\} \cup \big(\mathcal{I}_B \setminus \{P\}\big) \cup \big(M_B' \setminus \{vv^*\}\big) \cup U 
\end{array} 
\end{displaymath} 
%\noindent 
is an IPP of $G$. But then 
\begin{displaymath} 
\begin{array}{rrrcl} 
\text{ipp}(G) & \leq & |\mathcal{I}| & = & 1+(|\mathcal{I}_B|-1)+(|M_B'|-1)+|U| 
\\
\\
& & & = & |\mathcal{I}_B|+|M_B'|+|U|-1 
\\
\\
& & & \leq & |M_B|+|M_B'|+|U|-1 
\\
\\
& & & = & |M_B \cup M_B' \cup U|-1 
\\
\\
& & & = & |M \cup U|-1, 
\end{array} 
\end{displaymath} 
contrary to the fact that $M \cup U$ is a minimum IPP of $G$.~$\blacklozenge$ 

\medskip 

Since the block $B$ is odd (and also biconnected, by the definition of a block), Claim~\ref{prop-extremal-implies-one-even-block-claim-B-extremal} and Proposition~\ref{prop-odd-is-block} together imply that $B$ is an odd complete graph. Set $G' := G \setminus \big(V(B) \setminus \{v\})$. By Proposition~\ref{prop-delete-odd-leaf-clique}, $G'$ is IPP-extremal and the blocks of $G$ are precisely the block $B$ and the blocks of $G'$. Since $B$ is an odd complete graph, the result now follows immediately from the induction hypothesis applied to $G'$. 
\end{proof} 

\begin{proposition} \label{prop-G-extremal-implies-blocks-extremal} Let $G$ be a connected graph. Then both the following hold: 
\begin{enumerate}[(a)] 
\item if $G$ is IPP-extremal, then so are all its blocks; 
\item if $G$ is IndPP-extremal, then so are all its blocks. 
\end{enumerate} 
\end{proposition} 
\begin{proof} 
We may assume inductively that the proposition is true for graphs on fewer than $|V(G)|$ vertices. More precisely, we assume inductively that for all connected graphs $G'$ such that $|V(G')| < |V(G)|$, if $G'$ is IPP-extremal (resp.\ IndPP-extremal), then all blocks of $G'$ are IPP-extremal (resp.\ IndPP-extremal). We may further assume $G$ is not biconnected, for otherwise, $G$ has exactly one block (namely itself), and the result is immediate. Therefore, $G$ has at least two leaf-blocks. 

Suppose first that no leaf-block of $G$ is an odd complete graph. Since $G$ has at least two leaf-blocks, Proposition~\ref{prop-extremal-implies-one-even-block} now implies that $G$ is not IPP-extremal. Consequently, by Proposition~\ref{prop-indpp-implies-ipp}, $G$ is not IndPP-extremal either. Therefore, both (a) and (b) are vacuously true, and we are done. 

From now on, we assume that $G$ has a leaf-block, call it $B$, that is an odd complete graph. Let $v$ be the (unique) cut-vertex of $G$ that belongs to $B$. Set $G' := G \setminus \big(V(B) \setminus \{v\})$. Clearly, $G'$ is connected, and furthermore, Proposition~\ref{prop-delete-odd-leaf-clique} guarantees that the blocks of $G$ are precisely the block $B$ and the blocks of $G'$. Since $B$ is complete, it is clear that $B$ is both IPP-extremal and IndPP-extremal. On the other hand, Proposition~\ref{prop-delete-odd-leaf-clique} guarantees that if $G$ is IPP-extremal (resp.\ IndPP-extremal), then $G'$ is also IPP-extremal (resp.\ IndPP-extremal). The result now follows immediately from the induction hypothesis applied to $G'$. 
\end{proof} 

\medskip 

We are now ready to prove Theorems~\ref{thm-main} and~\ref{thm-main-indpp}. As we have previously pointed out, Proposition~\ref{prop-one-even-block-implies-extremal} yields the backward direction of Theorem~\ref{thm-main}, and Proposition~\ref{prop-one-even-block-implies-extremal-IndPP} yields the backward direction of Theorem~\ref{thm-main-indpp}. On the other hand, the forward direction of both theorems follows from Propositions~\ref{prop-extremal-implies-one-even-block} and~\ref{prop-G-extremal-implies-blocks-extremal}.\footnote{In the case of Theorem~\ref{thm-main-indpp}, we are also using the fact that IndPP-extremal graphs are IPP-extremal (by Proposition~\ref{prop-indpp-implies-ipp}), which means that Proposition~\ref{prop-extremal-implies-one-even-block} also applies to IndPP-extremal graphs.}

\section{Extremal graphs}\label{sec:corollaries}

In this section, we use Theorem~\ref{thm-main} (resp.\ Theorem~\ref{thm-main-indpp}) to obtain a complete characterization of connected odd graphs that are IPP-extremal (resp.\ IndPP-extremal), as well as a characterization of IPP-extremal block graphs (resp.\ IndPP-extremal block graphs). We begin with a simple proposition. 

\begin{proposition} \label{prop-even-block-even-card} If a connected graph $G$ has exactly $k$ even blocks, then $|V(G)|$ and $k+1$ have the same parity. In particular, any connected graph with no even blocks is odd, and any connected graph with exactly one even block is even. 
\end{proposition} 
\begin{proof} 
Fix a connected graph $G$ with exactly $k$ even blocks, and assume inductively that all connected graphs $G'$ such that $|V(G')| < |V(G)|$, if $G'$ has exactly $k'$ even blocks, then $|V(G')|$ and $k'+1$ have the same parity. We must show that $|V(G)|$ and $k+1$ have the same parity. 

If $G$ is biconnected, then $G$ has exactly one block, namely $G$ itself, and the result is immediate. So, we may assume that $G$ is not biconnected. Then $G$ has at least two leaf-blocks. Let $B$ be some leaf-block of $G$. Clearly, $|V(B)| \geq 2$. Let $v$ be the (unique) cut-vertex of $G$ that belongs to $B$, and set $G' := G \setminus \big(V(B) \setminus \{v\}\big)$. Then the blocks of $G$ are precisely the blocks of $G'$, plus the block $B$, and obviously, $G'$ is connected. Moreover, $|V(G)| = |V(G')|+|V(B)|-1$, and in particular, $|V(G')| < |V(G)|$. Let $k'$ be the number of even blocks of $G'$. By the induction hypothesis, $|V(G')|$ and $k'+1$ have the same parity. 

 Suppose first that the leaf-block $B$ is even. Then the number of even blocks of $G$ is $k = k'+1$, whereas $|V(G)|$ and $|V(G')|$ have the opposite parity (because $|V(G)| = |V(G')|+|V(B)|-1$, and $|V(B)|$ is even). Since $|V(G')|$ and $k'+1$ have the same parity, it follows that $|V(G)|$ and $k+1$ also have the same parity. 

 Suppose now that the leaf-block $B$ is odd. Then the number of even blocks of $G$ is $k = k'$, whereas $|V(G)|$ and $|V(G')|$ have the same parity (because $|V(G)| = |V(G')|+|V(B)|-1$, and $|V(B)|$ is odd). Since $|V(G')|$ and $k'+1$ have the same parity, so do $|V(G)|$ and $k+1$.
\end{proof}

Theorems~\ref{thm-main} and~\ref{thm-main-indpp}, together with  Proposition~\ref{prop-even-block-even-card}, readily yield the following four corollaries.

\begin{corollary} \label{cor-main-odd} Let $G$ be a connected odd graph. Then the following are equivalent: 
\begin{itemize} 
\item $G$ is IPP-extremal; 
\item $G$ is IndPP-extremal; 
\item every block of $G$ is an odd complete graph (and in particular, $G$ is a block graph). 
\end{itemize} 
\end{corollary} 
% \begin{proof} 
% If all blocks of $G$ are odd complete graphs, then Theorem~\ref{thm-main} guarantees that $\text{ipp}(G) = |V(G)|-\nu(G)$. Suppose, conversely, that $\text{ipp}(G) = |V(G)|-\nu(G)$. Since $G$ is connected and odd, Proposition~\ref{prop-even-block-even-card} guarantees that $G$ cannot have exactly one even block; so, by Theorem~\ref{thm-main}, all blocks of $G$ are odd complete graphs. 
% \end{proof} 

\begin{corollary} \label{cor-main-even} Let $G$ be a connected even graph. Then the following are equivalent: 
\begin{itemize} 
\item $G$ is IPP-extremal; 
\item $G$ has exactly one block that is not an odd complete graph, and that block is even and IPP-extremal. 
\end{itemize} 
\end{corollary} 
% \begin{proof} 
% Suppose first that $\text{ipp} = |V(G)|-\nu(G)$. Since $G$ is connected and even, Proposition~\ref{prop-even-block-even-card} guarantees that $G$ has at least one even block. But now Theorem~\ref{thm-main} guarantees that $G$ has exactly one block (call it $B$) that is not an odd complete graph, and the block $B$ is even and satisfies $\text{ipp}(B) = |V(B)|-\nu(B)$. The converse follows immediately from Theorem~\ref{thm-main}. 
% \end{proof} 

\begin{corollary} \label{cor-main-even-indpp} Let $G$ be a connected even graph. Then the following are equivalent: 
\begin{itemize} 
\item $G$ is IndPP-extremal; 
\item $G$ has exactly one block that is not an odd complete graph, and that block is even and IndPP-extremal. 
\end{itemize} 
\end{corollary} 

\begin{corollary} \label{cor-extremal-block-graph} Let $G$ be a connected block graph. Then the following are equivalent: 
\begin{itemize} 
\item $G$ is IPP-extremal; 
\item $G$ is IndPP-extremal; 
\item at most one block of $G$ is even. 
\end{itemize} 
\end{corollary} 
% \begin{proof} Since $G$ is a block graph, all its blocks are complete graphs. Moreover, it is clear that every complete graph $H$ satisfies $\text{ipp}(H) = |V(H)|-\nu(H)$. The result now follows immediately from Theorem~\ref{thm-main}. 
% \end{proof} 

\section{Concluding remarks and open questions}\label{sec:conclusion}

We have introduced an upper bound on the isometric path partition number and induced path partition number of a graph in terms of the matching number: any graph $G$ satisfies $\text{indpp}(G) \le \text{ipp}(G) \leq |V(G)|-\nu(G)$ (see Proposition~\ref{prop-upper-bound-trivial}). By Proposition~\ref{prop-components}, a graph is IPP-extremal (resp.\ IndPP-extremal) if and only if all its components are IPP-extremal (resp.\ IndPP-extremal). Theorems~\ref{thm-main} and~\ref{thm-main-indpp} provide a characterization of connected IPP-extremal graphs and IndPP-extremal graphs, respectively, in terms of their blocks. Corollary~\ref{cor-main-odd} gives a complete characterization of all connected odd graphs that are IPP-extremal and that are IndPP-extremal: a connected odd graph $G$ is IPP-extremal (resp.\ $G$ is IndPP-extremal) if and only if it is a block graph containing only odd blocks. On the other hand, by Corollary~\ref{cor-main-even} (resp.\ Corollary~\ref{cor-main-even-indpp}), a connected even graph $G$ is IPP-extremal (resp.\ IndPP-extremal) if and only if it contains exactly one block that is not a complete odd graph, and this one block is an even IPP-extremal graph (resp.\ IndPP-extremal graph). This reduces the problem of fully characterizing the IPP-extremal graphs and IndPP-extremal graphs to the even biconnected case. Thus, we propose the following questions. 

\begin{question} \label{q-biconn} 
Which biconnected even graphs are IPP-extremal?
\end{question}

\begin{question} \label{q-biconnIndPP} 
Which biconnected even graphs are IndPP-extremal?
\end{question}

We note that some obvious examples of biconnected even graphs that are IPP-extremal and IndPP-extremal include even complete graphs, the diamond, and the cycle $C_4$. 

Our proof of Theorems~\ref{thm-main} and~\ref{thm-main-indpp} made heavy use of inductive arguments. However, similar inductive arguments are unlikely to yield an answer to Questions~\ref{q-biconn} and~\ref{q-biconnIndPP}, and this is essentially because the IPP and IndPP numbers are not monotone with respect to the induced subgraph relation. More precisely, for an induced subgraph $H$ of a graph $G$, any one of the following is possible: $\text{ipp}(H) < \text{ipp}(G)$ (resp.\ $\text{indpp}(H) < \text{indpp}(G)$), $\text{ipp}(H) = \text{ipp}(G)$ (resp.\ $\text{indpp}(H) = \text{indpp}(G)$), $\text{ipp}(H) > \text{ipp}(G)$ (resp.\ $\text{indpp}(H) > \text{indpp}(G)$). It is easy to come up with examples of the first two possibilities (just consider complete graphs). For a concrete example of the third possibility, consider the graph $G$ and its induced subgraph $H$ represented in Figure~\ref{fig:example}. It is easy to see that $\text{indpp}(G) = \text{ipp}(G) = 2$, whereas $\text{indpp}(H) = \text{ipp}(H) = 3$. For the latter, we simply observe that no path of~$H$ can contain more than two leaves of~$H$; since $H$ contains five leaves, it follows that any IPP or IndPP of $H$ contains at least three paths, and consequently, $\text{ipp}(H) \geq 3$ and $\text{indpp}(H) \geq 3$. On the other hand, an IPP of $H$, which is also an IndPP of $H$, of cardinality $3$ is shown in Figure~\ref{fig:example}.

\begin{figure}[h]
\begin{center}
    \begin{tikzpicture}
    
 \node[minimum size=1.5cm, regular polygon, regular polygon sides=6, rotate=0] (hex) at (0,0) {};
        \foreach \x in {1,2,...,6}{%
            \node[mynode] at (hex.corner \x) (f\x) {};
    }
    \draw (f1)to(f2)to(f3)to(f4)to(f5)to(f6)to(f1);

\node[minimum size=3cm, regular polygon, regular polygon sides=6, rotate=0] (hex) at (0,0) {};

\foreach \x in {1,2,3,6}{%
            \node[mynode] at (hex.corner \x) (g\x) {};
    }
    \draw (f1)to(g1);
    \draw (f2)to(g2);
    \draw (f3)to(g3);
    \draw (f6)to(g6);
    %highlight
    \draw [-, line width=7pt, green, opacity=0.6] (g3) -- (f3) -- (f4) -- (f5) -- (f6) -- (g6);
    \draw [-, line width=7pt, yellow, opacity=0.5] (g2) -- (f2) -- (f1) -- (g1);
    
     \node at (0,-1.5) {$G$};

    \node[minimum size=1.5cm, regular polygon, regular polygon sides=6, rotate=0] (hex) at (5,0) {};
        \foreach \x in {1,2,3,4,6}{%
            \node[mynode] at (hex.corner \x) (f\x) {};
    }
    \draw (f1)to(f2)to(f3)to(f4);
    \draw (f6)to(f1);
    
\node[minimum size=3cm, regular polygon, regular polygon sides=6, rotate=0] (hex) at (5,0) {};

\foreach \x in {1,2,3,6}{%
            \node[mynode] at (hex.corner \x) (g\x) {};
    }
    \draw (f1)to(g1);
    \draw (f2)to(g2);
    \draw (f3)to(g3);
    \draw (f6)to(g6);
    %highlight
    
     \draw [-, line width=7pt, green, opacity=0.6] (g3) -- (f3) -- (f4);
    \draw [-, line width=7pt, yellow, opacity=0.5] (g2) -- (f2) -- (f1) -- (g1);
    \draw [-, line width=7pt, red, opacity=0.6] (g6) -- (f6);
    \node at (5,-1.5) {$H$};
        
    \end{tikzpicture}
    \end{center}
    \caption{An illustration of a graph $G$ and its subgraph $H$ where $\text{ipp}(G) < \text{ipp}(H)$ and $\text{indpp}(G) < \text{indpp}(H)$.} 
    \label{fig:example}
\end{figure}

%Other related parameters in this line of study include the (unrestricted) path partition (denoted as $\text{pp}$); which has been studied in~\cite{skupien1974path}. Notice that every isometric path is an induced path, and every induced path is a path; this, together with Proposition~\ref{prop-upper-bound-trivial} implies that, for any graph $G$, $\text{pp}(G) \leq \text{indpp}(G) \leq \text{ipp}(G) \leq |V(G)| - \nu(G)$. One can observe that every IndPP-extremal graph with respect to this upper bound is also IPP-extremal, but the converse is not true. A simple example for this case is a wheel graph on $6$-vertices (see Figure~\ref{fig:indpp}). On the other hand, the cycle $C_4$ is both IPP-extremal and IndPP-extremal, but not PP-extremal with respect to our upper bound. Furthermore, the only PP-extremal graphs are disjoint union of copies of $K_2$ and $K_1$; otherwise any connected component $|C| \geq 3$ allows a matching edge to be extended into a longer path, implying $pp(C) < |V(C)| - \nu(C)$. 

% Thus, we propose the following question. 

% \begin{question}
%     Find a characterization of graphs that are PP-extremal.
% \end{question}
Our upper bound for the IPP number raises an interesting algorithmic question as follows. 

\begin{question}
    Given a graph $G$ and an integer $k$, does there exist an IPP of $G$ of cardinality at most $|V(G)|-\nu(G)-k$?
\end{question}

Questions of this nature have been the subject of growing interest in the algorithmic community in recent years. Such formulations often form the foundation for parameterized complexity studies and are valuable in exploring the tractability and kernelizability of graph problems under certain structural constraints.

% In this context, it is worth pointing out that if $H$ is an induced subgraph of a graph $G$, then not every isometric path of $H$ need be\TSS{is, instead of need be} an isometric path of $G$, and in fact, it is possible that $\text{ipp}(H) > \text{ipp}(G)$. For a concrete example, take $G$ to be the graph obtained from a path $H$ on $k \geq 3$ vertices by adding a vertex $v$, adjacent to all vertices of $H$; then $\{H\}$ is an IPP of $H$ (because $H$ is an isometric path of itself), and consequently $\text{ipp}(H) = 1$, whereas $\text{diam}(G) = 2$, and so $\text{ipp}(G) \geq \Big\lceil \frac{|V(G)|}{3} \Big\rceil = \Big\lceil \frac{k+1}{3} \Big\rceil \geq 2$. 

\subsection*{Acknowledgment} 

We thank Prof. Florent Foucaud for the initial discussions on this problem. The original version of this manuscript dealt only with IPP. We thank
an anonymous referee for pointing out that most of our results about IPP generalize to IndPP, an observation that has been incorporated into this paper. We thank the other anonymous referee for pointing out to us that PP-extremal graphs are precisely the disjoint unions of $K_1$'s and $K_2$'s. We thank both referees for their careful reading of our manuscript and helpful suggestions.

\bibliographystyle{abbrv}
\bibliography{references-v1.1}

\section*{A list of corrections} 

Sections~2 and~3 of this manuscript are nearly identical to those of sections~2 and~3 of the published paper (cited on the first page of this manuscript). We have, however, corrected a few minor errors. Those that have mathematical significance are listed below. We do not list the purely linguistic and typographic changes that do not affect the mathematical content. 
\begin{enumerate} 
\item In the proof of Proposition~\ref{prop-components} (which corresponds to Proposition~1.5 of the published paper), $\mathcal{I}_i$ is a minimum IPP of $G_i$ (and not of $G$). 
\item In the second-to-last sentence of the proof of the Claim~\ref{prop-delete-edge-claim-delete-edge}, it should say ``and by construction, $xy$ belongs to this new matching'' rather than ``and by construction, $xy \in M$.'' 
\item For the sake of emphasis, in the second-to-last paragraph of the proof of Claim~\ref{prop-delete-edge-claim-xy-in-some-IPP-of-G}, we replaced ``of size $|\mathcal{I}|$'' by ``of size $|\mathcal{I}| = \text{ipp}(G)$.''
\item The proofs of Propositions~\ref{prop-one-even-block-implies-extremal} and~\ref{prop-one-even-block-implies-extremal-IndPP} should make no reference to Proposition~\ref{prop-odd-is-block}. Instead, the assertion in question follows simply from the assumptions of the two propositions (i.e.\ Propositions~\ref{prop-one-even-block-implies-extremal} and~\ref{prop-one-even-block-implies-extremal-IndPP}). In particular, the footnote from the proof of Proposition~\ref{prop-one-even-block-implies-extremal-IndPP} has been deleted. Moreover, at the end of the second-to-last paragraph of the proof of Proposition~\ref{prop-one-even-block-implies-extremal}, we added the following clarification: ``because clearly, all complete graphs are IPP-extremal.'' 
\item In the footnote from the proof of Claim~\ref{prop-v-extend-claim-NGv-cliques}, $C$ is supposed to be a connected component of $G[N_G(v)]$, and not of $G$. 
\item In the proof of Claim~\ref{prop-extremal-implies-one-even-block-claim-MB}, $B \setminus \{v,v^*\}$ is not necessarily an odd component of $G_U \setminus \{v,v^*\}$, but instead, some component of $B \setminus \{v,v^*\}$ is an odd component of $G_U \setminus \{v,v^*\}$. Moreover, $M \setminus \{vv^*\}$ is a perfect matching of $G_U \setminus \{v,v^*\}$, and not of $G_U$.
\item In the proof of Proposition~\ref{prop-v-extend} (in the paragraph right above Claim~\ref{prop-v-extend-claim-Q-spikes}), we replaced ``there exists an induced path in $Q := q_0q_1\dots q_tq_{t+1}$ ($t \geq 1$)'' by ``there exists an induced path $Q := q_0q_1\dots q_tq_{t+1}$ ($t \geq 1$) in $G \setminus v$.'' 
\item There was a typographic error in the proof of Claim~\ref{prop-extremal-implies-one-even-block-claim-B-extremal} of the published paper. It said ``and we deduce that is an IPP of $G$,'' whereas it should have stated: ``and we deduce that 
\begin{displaymath} 
\begin{array}{rcl} 
\mathcal{I} & := & \{v^*vp_1\dots p_t\} \cup \big(\mathcal{I}_B \setminus \{P\}\big) \cup \big(M_B' \setminus \{vv^*\}\big) \cup U
\end{array} 
\end{displaymath} 
is an IPP of $G$.'' (In the published paper, the definition of $\mathcal{I}$ was moved to an incorrect place on the page.) 
Further, the displayed chain of inequalities should begin with
$\text{ipp}(G) \leq |\mathcal{I}|$, rather than
$\text{ipp}(B) \leq |\mathcal{I}|$, since $\mathcal{I}$ is an IPP of $G$
and $M \cup U$ is a minimum IPP of $G$.
\item In the third paragraph of the proof of Proposition~\ref{prop-G-extremal-implies-blocks-extremal}, the following sentence was added as the second sentence: ``Let $v$ be the (unique) cut-vertex of $G$ that belongs to $B$.''
\end{enumerate} 

\end{document}